\documentclass[11pt]{amsart}

\usepackage{amsfonts,amsmath,amsthm,amscd,amssymb,latexsym,mathtools,amsxtra} 

\usepackage{textcomp} 

\usepackage{mathrsfs}
\usepackage{txfonts}
\usepackage{bbm}
\usepackage{wasysym}
\usepackage{psfrag}
\usepackage{graphics,epsfig}  
\usepackage{comment}
\usepackage{enumerate}
\usepackage{url}
\usepackage{color}
\usepackage[FIGTOPCAP]{subfigure}

\input{macros_MATHmarta.sty}

\newtheorem{theorem}{Theorem}[section]

\newtheorem{definition}[theorem]{Definition}

\newtheorem{remark}[theorem]{Remark}

\RequirePackage{fix-cm}
\graphicspath{}

\begin{document}


\title[non-twist tori]{non-twist tori in conformally symplectic systems}

\author[R. Calleja]{Renato Calleja$^{*}$}
\address{XX} 
\email{calleja@mym.iimas.unam.mx}

\author[M. Canadell]{Marta Canadell$^{**}$}
\address{Institute for Computational and Experimental Research in Mathematics (ICERM),
 Brown University. 121 South Main St., Providence (RI), USA} 
\email{marta\_canadell@brown.edu} 

\author[A. Haro]{Alex Haro$^{***}$}
\address{Departament de Matem\`atiques i Inform\`atica, Universitat de Barcelona,
 Gran Via 585, 08007 Barcelona, Spain} 
\address{Barcelona Graduate School of Mathematics (BGSMath)}
\email{alex@maia.ub.es}

\date{\today}


\thanks{M.C. and A.H. partially supported by MTM2012-32541 and 2014-SGR-1145. 
M.C. also supported by NSF DMS-1162544 and XX. 
R.C was partially supported by DGAPA-UNAM projects PAPIIT
IA 102818, IN101020 and by UIU project UCM-04-2019.\\
\indent\emph{E-mail addresses:} 
$^{*} $\texttt{calleja@mym.iimas.unam.mx}, 
$^{**}$\texttt{marta\_canadell@brown.edu},
$^{***}$\texttt{alex@maia.ub.es}
}

\maketitle


\begin{abstract}
  Dissipative mechanical systems on the torus with a friction that is proportional
  to the velocity are modeled 
  by conformally symplectic maps on the annulus, which are maps that transport the symplectic form 
  into a multiple of itself (with a conformal factor smaller than $1$). It is important to understand
  the structure and the 
  dynamics on the attractors. With the aid of parameters, and under suitable non-degeneracy conditions, 
  one can obtain that, by adjusting parameters, there is an attractor that is an
  invariant torus whose internal dynamics is conjugate to a rotation \cite{CallejaCL13a}.  By analogy
  with symplectic dynamics, there have been some debate 
  in establishing appropriate  definitions for twist and non-twist invariant tori (or systems). The
  purpose of this paper is 
  two-fold: (a) to establish proper definitions of twist and non-twist invariant tori in families of
  conformally symplectic systems;
  (b) to derive algorithms of computation of non-twist invariant tori. 
  The last part of the paper is devoted to implementations of the algorithms, illustrating the definitions presented in 
  this paper, and exploring the mechanisms of breakdown of non-twist tori.
  For the sake of simplicity we have considered here 2D systems, i.e. defined in the 2D annulus, but generalization 
  to higher dimensions is straightforward. 
\end{abstract}


\def\AUE{{A\times U\times E}}
\def\afun{\mbox{ $\vartheta$}}

\section{The Introduction}

Conformally symplectic systems model some mechanical systems with dissipation, 
in which the friction is proportional to the velocity.
Geometrically, conformally symplectic systems transport a symplectic form into a
multiple of itself. When the conformal factor is less than one the systems
contract the form and are dissipative.
In contrast to symplectic systems, dissipative systems have attractors.
Although dissipative systems have less asymptotic behaviors by themselves, one recovers
asymptotic behaviors by adding adjusting parameters.
There has been a lot of interest in the case these attractors are invariant smooth
tori 
that contain quasi-periodic dynamics. Obtaining quasi-periodic dynamics is proved thanks to
the presence of parameters in the system
and some non-degeneracy condition that is referred to as twist condition in analogy 
of the common twist condition that appears in symplectic systems, \cite{Moser66b, Moser67,
  BroerHS96, CallejaCL13a, CanadellH17a}.
In this paper we are interested in developing algorithms for computing quasi-periodic
circles when an analogue of the twist condition fails. In fact, a first task is to identify
the proper definition for a non-twist circle in this context. 

To the best of our knowledge, this paper presents a first attempt for considering
non-twist tori in dissipative systems. We
will present algorithms for 
the simplest 2D case. We have not proved here the convergence of the algorithms, but just
applied them in several examples. However, we expect that a proof could be done by using
standard KAM techniques, see for example \cite{BroerHS96, GonzalezJLV05, CallejaCL13a, GonzalezHL13,
  CanadellH17a}.

\paragraph*{Organization of the paper} In Section~\ref{section:setting} we introduce the setting, and present
suitable definitions of  non-twist tori in the context of conformally symplectic dynamics. In Section~\ref{section:algorithms} we 
describe a methodology for the computation of invariant tori in conformally symplectic systems, and, more importantly, 
for the computation and continuation of non-twist invariant tori in these systems. Section~\ref{section:applications} is devoted 
to implementations of the algorithms to several examples, referred to as dissipative standard non-twist families,
illustrating the concepts and algorithms presented in 
this paper, and to the analysis of the breakdown of non-twist invariant tori.


\section{The definitions}
\label{section:setting}

In this section, we present and motivate the definition of non-twist tori
in families of conformally symplectic
systems in a 2-dimensional phase space.
The phase space is the annulus ${\TT\times\RR}$, 
endowed with coordinates $z= (x,y)$, being $\TT= \RR/\ZZ$ the torus. 
We consider a $3$-parameter family of (dissipative) conformally symplectic maps in the annulus 
$\TT\times\RR$, with conformal factor $\sigma\in ]0,1[$, 
given by a smooth map 
$F= (F^x,F^y):(\TT\times\RR)\times \AUE\to\TT\times\RR$, where $A,U,E\subset\RR$ are open intervals,
for which we will write 
\begin{equation*}
F_{a,\mu,\eps}(z):=F(z;a,\mu,\eps),
\end{equation*}
such that, for each $(a,\mu,\eps) \in \AUE$:
\begin{itemize}
\item $F_{a,\mu,\eps}$ is a diffeomorphism, homotopic to the identity (i.e., lifting to the covering space, 
$F^x(x,y)-x$ is a periodic function);
\item for all $z\in \TT\times\RR$, $\det {\rm D}F_{a,\mu,\eps}(z)= \sigma$.
\end{itemize}
Notice that, by considering the symplectic product given by matrix 
\[
\Omega=
\begin{pmatrix}
0 & -1\\
1 & 0
\end{pmatrix}, 
\]
the determinant condition may be written as 
\[
	({\rm D}F_{a,\mu,\eps}(z))^\top\ \Omega \ {\rm D}F_{a,\mu,\eps}(z) = \sigma \Omega.
\]
In the limiting case $\sigma= 1$, these diffeomorphisms are symplectic.
In the sequel, the parameters $a,\mu,\eps$ will play  very different roles. 

We are interested in the continuation of normally (hyperbolic) contracting invariant rotational circles for $F$ 
with respect to parameters $a,\mu,\eps$, and most particularly, those for which the internal dynamics is quasi-periodic 
(with a certain fixed Diophantine rotation number $\omega$) and are degenerate in some sense to be specified below. The non-degenerate case has been considered in \cite{CallejaC10,CallejaCL13a,CallejaF12}. 

\begin{remark}
Conformal sympleticity imposes severe restrictions for the existence of invariant circles. For instance, there can not exist 
invariant librational circles (those that are homotopically trivial), and no more than one invariant rotational circle.
\end{remark}

\begin{remark}
The case in which the conformal factor $\sigma(z;a,\mu,\eps)= \det {\rm D}F_{a,\mu,\eps}(z)$ 
depends on $z$ and the parameters  could be also considered with 
slight modifications of the arguments. Here, we consider the constant case for the sake of simplicity.
\end{remark}

\begin{definition}
 \label{def:InvKf}
Given fixed parameters $(a_0,\mu_0,\eps_0)\in\AUE$, we say that the circle $\calK$ parameterized by 
$K:\TT\to \TT\times\RR$ is a \emph{$F_{a_0,\mu_0,\eps_0}$-invariant rotational circle} with 
internal dynamics $f:\TT\to\TT$, if $\calK$ is homotopic to the zero section (so that, 
$\theta\to K^x(\theta)-\theta$ is 1-periodic) and
the couple $(K,f)$ satisfies the invariance equation: 
\begin{equation}
F_{a_0,\mu_0,\eps_0}( K(\theta)) - K(f(\theta)) = 0.
\label{eq:InvKf}
\end{equation}
Notice that the internal dynamics is homotopic to the identity (hence, the lift of $f(\theta)-\theta$ is $1$-periodic).

We say that, moreover, $\calK$ is normally contracting if there exists positive constants $C,\lambda<1$ such that, 
for all $\theta\in \TT$, $k\in \ZZ^+$, 
\begin{equation}
\label{def:contraction}
	\sigma^k \left(f'(f^{k-1}(\theta)) \dots f'(f(\theta)) f'(\theta)\right)^{-2} \leq C \lambda^k,
\end{equation}
where $'$ denotes the derivative with respect to $\theta$.

A particular case, is when the internal dynamics is (smoothly) conjugate to a rotation by a certain angle $\omega\in \RR\setminus\QQ$ and, hence, 
we can reparametize $\calK$ so that 
\begin{equation}
F_{a_0,\mu_0,\eps_0}( K(\theta)) - K(\theta+\omega) = 0,
\label{eq:InvKq}
\end{equation}
that is $f(\theta)= \theta+\omega$.
We will say then that   $\calK$ is a \emph{quasi-periodic $F_{a_0,\mu_0,\eps_0}$-invariant rotational circle}. 
\end{definition}

\begin{remark}
Notice that, by a change of variables, we can assume the phase condition 
\begin{equation}\label{eq:Phasecondi}
\avg{K^x(\theta)-\theta} =  0.
\end{equation}
More especifically, for a reparameterization $K_\alpha(\theta)= K(\theta+\alpha)$, for which the corresponding internal dynamics is given by $f_\alpha(\theta)= f(\theta+\alpha)-\alpha$, 
the phase is $\avg{K^x_\alpha(\theta)-\theta} = \alpha + \avg{K^x(\theta)-\theta}$, and hence 
one can adjust $\alpha$ in order to adjust the phase condition.
\end{remark}

The meaning of the normally contracting property is that there is a normal invariant bundle for which the linearized dynamics is contracting, and whose rate of contraction dominates the internal dynamics on the circle. In the sequel, we will formulate this idea in a rather computational way, since 
we are interested here in numerical algorithms and their implementation. 
The tangent bundle ${\mathcal T}{\calK}$ to the circle $\calK$ is spanned by the derivative map $K':\TT\to \RR\times\RR$ . 
We can consider a normal bundle ${\mathcal N}_0{\calK}$ generated by $N_0:\TT\to \RR\times\RR$, where 
\begin{equation}
\label{eq:defN0}
	N_0(\theta)= \Omega\ K'(\theta)\ (K'(\theta)^\top K'(\theta))^{-1}.
\end{equation}
Notice that, with this choice
\[
     N_0(\theta)^\top \Omega\ K'(\theta)= 1.
\]
The geometrical meaning is that the area of the parallelogram generated by $K'(\theta)$ and $N_0(\theta)$ is $1$.
While the tangent bundle is invariant for the linearized dynamics,  and in particular 
\[
     \Dif F_{a_0,\mu_0,\eps_0}(K(\theta)) K'(\theta) = K'(f(\theta)) f'(\theta),
\]
the normal bundle ${\mathcal N}_0{\calK}$ could be non-invariant, since 
\[
     \Dif F_{a_0,\mu_0,\eps_0}(K(\theta)) N_0(\theta) = K'(f(\theta))\ t_0(\theta) + N_0(f(\theta))\ \frac{\sigma}{f'(\theta)},
\]
where
\begin{equation}
\label{eq:deft0}
   t_0(\theta) = N_0(f(\theta))^\top \Omega \ \Dif F_{a_0,\mu_0,\eps_0}(K(\theta)) N_0(\theta).
\end{equation}
In order to construct an invariant normal bundle ${\mathcal N}{\calK}$ spanned by a suitable $N:\TT\to \RR\times\RR$, we write
\begin{equation}
N(\theta) = K'(\theta) \afun(\theta) + N_0(\theta),
\label{eq:defN}
\end{equation}
for which $N(\theta)^\top \Omega\ K'(\theta)= 1$, 
and realize that 
\[
     \Dif F_{a_0,\mu_0,\eps_0}(K(\theta)) N(\theta) = K'(f(\theta))\ t(\theta) + N(f(\theta))\ \frac{\sigma}{f'(\theta)},
\]
where
\[
   t(\theta) = t_0(\theta) + f'(\theta) \afun(\theta) - \frac{\sigma}{f'(\theta)} \afun(f(\theta)).
\]
Hence, we  make $t(\theta)= 0$ 
by taking 
\[
\afun(\theta)= -\sum_{k= 0}^\infty \frac{\sigma^k}{\left(f'(f^{k-1}(\theta)) \dots f'(\theta)\right)^2} \cdot \frac{t_0(f^k(\theta))}{f'(f^k(\theta))}.
\]
Notice the convergence of the series is guaranteed by the normal contraction condition \eqref{def:contraction}.

In a nutshell, we have just constructed a frame $P:\TT\to\RR^{2\times 2}$, defined by juxtaposing $K'$ and $N$,  i.e. 
\begin{equation}
  \label{adapted_frame}
  P(\theta)= \begin{pmatrix} K'(\theta) & N(\theta) \end{pmatrix},
\end{equation}
that satisfies $\det P(\theta)= 1$ (the frame is symplectic) and reduces the linearized dynamics to diagonal form:
\[
	P(f(\theta))^{-1} \Dif F_{a_0,\mu_0,\eps_0}(K(\theta)) P(\theta) = 
	\begin{pmatrix} f'(\theta) & 0 \\ 0 & \tfrac{\sigma}{f'(\theta)} \end{pmatrix}.
\]

\begin{remark}
The rate of contraction $\lambda$ is a dynamical observable of the contracting condition, and has to be $\lambda<1$.
Another important observable that measures the quality of the hyperbolicity 
property is the (minimum) angle between the invariant bundles. In the setting of the present paper, this is given by 
  \begin{equation}
      \alpha = \min_{\theta\in\TT} \left|\arctan \left( \frac{1}{\vartheta(\theta) [K'(\theta)^\top K'(\theta)]} \right)\right|. 
  \end{equation}
 In the case $\alpha > 0$, there is a well-defined splitting in tangent and invariant normal bundles.
\end{remark}

\begin{remark}
We will be mainly interested in the quasi-periodic case, for which $f(\theta)= \theta+\omega$ (see \eqref{eq:InvKq}). 
Hence, $f'(\theta)= 1$, and the rate of contraction is $\lambda= \sigma<1$. In this case, the quality of the hyperbolicity 
property is essentially given by the positiveness of $\alpha$. 
\end{remark}

It is well-known that normal contractiveness (and, more generally, normal hyperbolicity is an open condition. 
Hence, if ${\calK}_{a_0,\mu_0,\eps_0}$ 
is a normally contracting $F_{a_0,\mu_0,\eps_0}$-invariant rotational circle, then there is an open neighborhood of 
$(a_0,\mu_0,\eps_0)$ in $\AUE$ 
for which there is a normally contracting $F_{a,\mu,\eps}$-invariant rotational circle ${\calK}_{a,\mu,\eps}$ for each $(a,\mu,\eps)$ in such a neighborhood. Without loss of generality, we consider this neighborhood to be also $\AUE$. Hence, we assume there are smooth maps $K:\TT\times\AUE \to \TT\to\RR$ and $f:\TT\times\AUE\to\TT$ such 
that, for each $(a,\mu,\eps)\in \AUE$:
\[
	F_{a,\mu,\eps}(K_{a,\mu,\eps}(\theta)) - K_{a,\mu,\eps}(f_{a,\mu,\eps}(\theta))= 0.
\]	 
Hence, we can define a rotation number function $\rho:\AUE\to \RR$ such that, for each 
$(a,\mu,\eps)\in \AUE$, $\rho(a,\mu,\eps)$ is 
the rotation number of $f_{a,\mu,\eps}$. The resonant set is the set of parameters for which the rotation number is rational
(and the internal dynamics possesses periodic orbits), 
and the non-resonant set corresponds to irrational rotation numbers. The regularity 
of the rotational circles jump from being (typically) finitely differentiable in the resonant set (generically with non-empty interior) to $C^\infty$ 
(or even real-analytic) in the non-resonant set (generically with empty interior).

In this setting, we introduce twist (and non-twist) notions for quasi-periodic invariant rotational circles, 
and emphasize the different roles of parameters. 
For $(\mu_0,\eps_0)\in U\times E$ fixed, we have a one-parameter family of  invariant rotational circles 
$a\to K_{a,\mu_0,\eps_0}$, and the corresponding rotation number function $a\to \rho(a,\mu_0,\eps_0)$ is (typically) a devil staircase, whose steps correspond to resonances.
Then, we say that an invariant rotational circle $K_{a_0,\mu_0,\eps_0}$ that is quasi-periodic satisfies the twist condition  with respect to parameter $a$ (or that it is an $a$-twist invariant rotational circle)
if the devil staircase is  strictly monotone  at $a_0$, otherwise we say it is non-twist with
respect to parameter $a$ 
(or that it is an non-$a$-twist invariant rotational circle). In place of this 
dynamical definition of the twist property we will use the following analytical definition, which is more practical.

\begin{definition}
 \label{def:InvK_NT}
Let $\calK$ be a quasi-periodic $F_{a_0,\mu_0,\eps_0}$-invariant rotational circle $\calK$ parameterized by $K:\TT\to \TT\times\RR$, with 
irrational rotation number $\omega$, so that 
\begin{equation}
\label{eq:invariance}
F_{a_0,\mu_0,\eps_0}(K(\theta))-K(\theta+\omega) = 0.
\end{equation}
We define the $a$-twist (the twist with respect to parameter $a$) 
to be the number 
\begin{equation}
\label{eq:a-twist}
	b_a(K; a_0,\mu_0,\eps_0)= \avg{ N(\theta+\omega)^\top \Omega\ \Dif_a F_{a_0,\mu_0,\eps_0}(K(\theta))},
\end{equation}
where $N$ is the parameterization of the normal invariant bundle.
Then, the circle is $a$-twist if $b_a(K; a_0,\mu_0,\eps_0)\neq 0$, and 
non-$a$-twist if $b_a(K; a_0,\mu_0,\eps_0)= 0$.
\end{definition}

In the previous definition, the twist property has to do with the fact that, for fixed
$(\mu_0,\eps_0)$, 
one can tune parameter $a$ to fix the 
rotation number to be $\omega$ (since the rotation number is a strictly monotone function of $a$
around $a_0$). In fact, for $\omega$ Diophantine, KAM techniques can be applied to obtain
real-analytic solutions of the invariance equation, under suitable sufficient conditions
including the twist condition, as demonstrated in \cite{CallejaCL13a}. Hence, by using the
  implicit function theorem, one obtains a mapping $a= a(\mu,\eps)$ such that the invariant
rotational circles  ${\calK}_{a(\mu,\eps),\mu,\eps}$ are real-analytic and their dynamics is a rotation
by $\omega$. In this context, we refer to $a$ as the {\em adjusting parameter}. We note that
there are situations in which the adjusting parameter $a$
does not change (to first order) the rotation number 
since the $a$-twist condition (with respect to such parameter) fails. 

In this paper we are interested in studying the boundaries of twist property, particularly in
developing algorithms for computing rotational invariant
circles with internal dynamics given by the rotation by $\omega$
when the twist condition, with respect to a parameter, fails.  
Notice that, besides parameter $a$ that is in principle designed to adjust the dynamics to be
the rotation by $\omega$, we 
need an extra parameter to select degenerate (with respect to $a$) invariant
rotational circles. This is the role 
of the {\em unfolding parameter $\mu$}. In some sense, $\mu$ is the truly adjusting parameter,
since $a$ can not do the job, and, given $a$ and $\eps$,
we can select a $\mu= \mu(a,\eps)$ and an invariant circle parameterized by a certain
$K_{a, \mu(a,\eps),\eps}$ so that 
the internal dynamics is a rigid rotation $f_{a,\mu(a,\eps),\eps}(\theta)= \theta+\omega$.
By writing $\bar b_a(a,\eps)= b_a(K_{a, \mu(a,\eps),\eps}, a, \mu(a,\eps),\eps)$ we aim to solve
the equation
\[
 \bar b_a(a,\eps_0)= 0, 
\]
for each $\eps_0$ fixed. Hence, $\eps$ is a continuation parameter and the goal
is finding $a= a(\eps)$ 
(and then $\mu= \mu(a(\eps),\eps)$) so that 
$\bar b_a(a(\eps),\eps)= 0$, in order to obtain a curve of
non-$a$-twist circles in the parameter space 
$\AUE$. By the  implicit function theorem, a sufficient condition is that 
$\bar b_a(a_0,\eps_0)= 0$ and $\frac{\partial \bar b_a}{\partial a}(a_0,\eps_0) \neq 0$
for a given $(a_0,\eps_0)$.

\section{The algorithms}
\label{section:algorithms}

In this section we present several algorithms for computing invariant rotational circles of 2-dimensional 
conformally symplectic systems. Firstly, tailoring the algorithms presented in \cite{Canadell14_thesis,HaroCFLM16} 
(see also \cite{Granados17} for other implementations), we give an algorithm to compute  normally contracting
circles, with non-fixed dynamics and hence, no parameters are needed. This will be useful to illustrate the
behavior of the rotation number function and its link to the twist
property for systems depending on parameters. Algorithms for computing invariant circles with fixed
quasi-periodic dynamics under twist conditions with respect to parameters (in particular, with respect
to the adjusting parameter $a$) are presented in \cite{CallejaC10,CallejaF12,CanadellH14,CanadellH17b}. 
We present here algorithms for computing non-$a$-twist invariant rotational circles
(adding the unfolding parameter $\mu$), and continuation (with respect to perturbing parameter $\eps$).
The  algorithms to solve the invariance equations are based on Newton's method.

\subsection{General algorithm for invariant tori in conformally symplectic systems}
\label{sec:GeneralConf}

In this section, we give an algorithm to compute the parameterization of an invariant rotational circle 
$\calK$ parameterized by $K:\TT\to\TT\times\RR$, for a conformally symplectic diffeomorphism 
$F:\TT\times\RR\to\TT\times\RR$,
with its unknown internal dynamics $f:\TT\to\TT$ (we emphasize the absence of parameters, since no tuning is needed). 
In particular, we explain how to perform one step of a Newton-like method to solve the invariance equation~\eqref{eq:InvKf}.

Let us assume that $K$ is approximately invariant, and let  $E:\TT\to \RR^2$ be the invariance error function given by
\begin{equation}
E(\theta)= F(K(\theta))-K(f(\theta)).\\
\label{eq:FullSystem_f}
\end{equation}
The goal of one step of Newton's method is to compute the corrections 
$\Delta K:\TT\to\RR^2,\Delta f:\TT\to\RR$ of $K$ and $f$,
respectively. The functions $\Delta K$ and $\Delta f$ are
$1$-periodic functions that are given in such a way that the error
estimates of the new approximations 
$\bar K= K + \Delta K, \bar f= f + \Delta f$,
are quadratically smaller with respect to the initial error
estimates. Then, by substituting the new approximations 
of the invariance equation \eqref{eq:InvKf}, using first order
Taylor expansion, we obtain:
\begin{equation*}
\begin{split}
  0  & = F(K(\theta)+\Delta K(\theta))- K(f(\theta)+\Deltaf(\theta)) - \Delta K(f(\theta)+\Deltaf(\theta))\\
     & = E(\theta)+\Dif F(K(\theta))\Delta K(\theta) -\Dif K (f(\theta))\Deltaf(\theta) - \DeltaK(f(\theta)) +\mathcal{O}_2, 
    \end{split}
\end{equation*}
where 
$\mathcal{O}_2$
includes the second order terms.  Hence, 
in principle, the Newton step consists in solving the linearized equation 
\begin{equation}
\label{Newton1}
\begin{split}
  \Dif F(K(\theta))\Delta K(\theta) -\Dif K (f(\theta))\Deltaf(\theta) - \DeltaK(f(\theta)) &= -E(\theta). 
    \end{split}
\end{equation}
Instead, we will solve the previous equation with an error that is
quadratically small with respect to 
the invariance error, $E$. To do so, we use adapted frames as follows.


First, we compute the frame $P(\theta)$ given in \eqref{adapted_frame}. Since the circle $\calK$ is not a priori invariant, 
then the linearized dynamics is approximately reduced to the diagonal form
\[
	\Lambda(\theta)= \begin{pmatrix} f'(\theta) & 0 \\ 0 & \tfrac{\sigma}{f'(\theta)} \end{pmatrix}.
\]
More especifically, there is a reducibility error function $E_r:\TT\to \RR^{2\times 2}$, given by
\begin{equation}
\label{approximate_reducibility}
E_r(\theta)= {\Dif}F(K(\theta)) P(\theta)- P(f(\theta)) \Lambda(\theta).
\end{equation}
One obtains
\[
E_r(\theta) = \begin{pmatrix} E'(\theta) & E_r^N(\theta) \end{pmatrix}
\]
where
\[
	E_r^N(\theta)= \tfrac{1}{f'(\theta)} \left(E'(\theta)^\top \Omega \Dif F(K(\theta)) N_0(\theta) \right) N_0(f(\theta)) 
	+ E'(\theta) \afun(\theta),
\]
being $N_0(\theta)$ given by \eqref{eq:defN0}.

Second, we write the correction term of the parameterization of the rotational circle as 
$ \Delta K(\theta)=  P(\theta)\xi(\theta)$,
where $\xi:\TT\to \RR^2$ is a periodic function. 
Then, by multiplying \eqref{Newton1}  by $P(f(\theta))^{-1}$,  using approximate reducibility \eqref{approximate_reducibility} 
and neglecting quadratically small terms, we obtain the following cohomological equation 
\begin{equation}
\label{eq:Newton_general}
\Lambda(\theta)\xi(\theta)-\xi(f(\theta))
-\left(\begin{array}{c}
\Deltaf(\theta)\\ 0
\end{array}\right) = \eta(\theta),
\end{equation}
where $\eta(\theta)= -P(f(\theta))^{-1} E(\theta) $ is the error of invariance  in the adapted frame. 

Third, since $\Lambda$ is diagonal, we split Equation~\eqref{eq:Newton_general} into tangent and normal (stable)
components, so we obtain the following two cohomological equations:
\begin{eqnarray}
f'(\theta)\xiL(\theta)-\xiL(f(\theta))-\Deltaf(\theta) &=& \etaL(\theta),
\label{eq:Cohom_tan} \\
\frac{\sigma}{f'(\theta)}\xiN(\theta)-\xiN(f(\theta)) &=& \etaN(\theta).
\label{eq:Cohom_stab}
\end{eqnarray}
Hence, we solve Equation~\eqref{eq:Cohom_stab} by simple iteration for the fixed point equation
\begin{equation*}
\xiN(\theta) =-\etaN(f^{-1}(\theta)) + \frac{\sigma}{f'(f^{-1}(\theta))}\xiN(f^{-1}(\theta)).
\end{equation*}
On the other hand, to solve Equation~\eqref{eq:Cohom_tan}, we need to solve  an \emph{underdetermined equation}: one equation with two unknowns ($\xiL$ and $\Deltaf$ ). 
Then, the simplest choice to solve this Equation~\eqref{eq:Cohom_tan} is by choosing the solution given by
\begin{equation*}
\xiL(\theta)=0,\  \Deltaf(\theta)=-\etaL(\theta).
\end{equation*}

Four, and last, we obtain the new approximations 
\[
\bar K(\theta)= K(\theta) +  N(\theta) \xiN(\theta),\  \bar f(\theta)= f (\theta)+ \Deltaf(\theta).
\]
With this computation we finish one step of Newton method.

The implementation of Newton method in a computer starts by choosing a method of representation of periodic functions. There are several methods at hand, such as trigonometric polynomials (via FFT), splines or (local) interpolating polynomials. While trigonometric polynomials are especially adapted to rigid rotations, and they will be used later when fixing quasi-periodic dynamics, in the present case they are computationally expensive and we  have used (local) interpolating polynomials \cite{Canadell14_thesis,HaroCFLM16} (for use of splines, see \cite{Granados17}).
We emphasize that when implementing the continuation with respect to parameters of the invariant rotational circles 
using derivatives with respect to such parameter one has to face linearized equations of the same type we have explained
in this section.


\bigskip \bigskip
\subsection{Algorithm for non-twist invariant tori in conformally symplectic systems}
\label{sec:NT-tori}
In this section, we give an algorithm to compute non-twist (with respect to a parameter)
invariant rotational circles with fixed frequency $\omega$, which is assumed to be
Diophantine. In fact, we present an algorithm to fix the $a$-twist $b_a$
  (see \eqref{eq:a-twist})
to a given value $b_a^0$ (with $b_a^0= 0$ in the non-$a$-twist case).
In contrast with the general algorithm described in Section~\ref{sec:GeneralConf},
we do not perform corrections of the internal dynamics, but adjust parameters
$a$ and $\mu$ (for each $\eps$) to fix the dynamics to the given rotation
and the $a$-twist. 
Hence, in the context of the introduction, 
we fix perturbation parameter $\eps_0$ and look for solutions $(K,a,\mu)$ of the system of equations
\begin{eqnarray}
\label{eq:Invariance0}
F(K(\theta),a,\mu,\eps_0)-K(\theta+\omega) &=& 0,\\
\label{eq:Phase0}
\avg{K^x(\theta)-\theta} & = & 0,\\
\label{eq:Twist0}
b_a(K;a,\mu,\eps_0) - b_a^0 & = & 0.
\end{eqnarray}

Assume then we have an approximate solution $(K,a,\mu)$  of \eqref{eq:Invariance0}, \eqref{eq:Phase0}, \eqref{eq:Twist0}.  
As we emphasize in the previous section, the aim to perform one step of the Newton's method is computing the corrections 
$(\Delta K, \Delta a, \Delta \mu)$ to obtain a new approximate solution $(\bar K,\bar a,\bar \mu)$ which will have an error that is quadratically small with respect to the initial error, even though the linearized equations are solved approximately using appropriate frames. Hence, the starting point is a triple $(K,a,\mu)$ such that
\begin{eqnarray}
\label{eq:InvarianceError}
F(K(\theta),a,\mu,\eps_0)-K(\theta+\omega) &=& E(\theta),\\
\label{eq:PhaseError}
\avg{K^x(\theta)-\theta} & = & e_p,\\
\label{eq:TwistError}
b_a(K;a,\mu,\eps_0) - b_a^0 & = & e_b,
\end{eqnarray}
where $E:\TT\to\RR^2$ and $e_p, e_b$ are small enough. In the following, we will proceed in two steps: 1) for any 
$\Delta a$, we compute $\Delta K$ and $\Delta \mu$ to improve \eqref{eq:InvarianceError} and \eqref{eq:PhaseError};
2) we adjust $\Delta a$ (and hence $\Delta K$ and $\Delta \mu$) to improve \eqref{eq:TwistError}.

The first part consists essentially in solving approximately 
\begin{equation}
\label{eq:Newton_qp}
\begin{split}
\Dif F(K(\theta),a,\mu,\eps_0)\Delta K(\theta)
      +\Dif_a F(K(\theta),a,\mu,\eps_0)\Delta a +\Dif_\mu F(K(\theta),a,\mu,\eps_0)\Delta \mu
      - \DeltaK(\theta+\omega)  &=  -E(\theta), 
\end{split}
\end{equation}
for any $\Delta a$ (adjusting also the phase condition, although this step could be done at the end of the iteration).
This is again performed with the aid of an adapted frame.

To do so, we first compute, from $L(\theta)= K'(\theta)$,  the expressions of $N_0(\theta)$ and $t_0(\theta)$, $\afun(\theta)$ and 
$N(\theta)$, and then the frame $P:\TT\to\RR^{2\times 2}$ given by
\[
P(\theta)= \begin{pmatrix} K'(\theta) & N(\theta) \end{pmatrix},
\]
that satisfies $\det P(\theta)= 1$ and 
\begin{equation}
\label{approximate_reducibility_qp}
	{\Dif}F(K(\theta),a,\mu,\eps_0) P(\theta)= P(\theta+\omega) \Lambda(\theta) + E_r(\theta),
\end{equation}
where 
\[
	\Lambda(\theta)= \begin{pmatrix} 1 & 0 \\ 0 & \sigma \end{pmatrix}
\]
and the reducibility error is 
\[
	E_r(\theta)= \begin{pmatrix} E'(\theta) & E_r^N(\theta) \end{pmatrix},
\]
with 
\[
	E_r^N(\theta)= \left(E'(\theta)^\top \Omega \ \Dif F(K(\theta),a,\mu,\eps_0) N_0(\theta) \right) N_0(\theta+\omega) 
	+ E'(\theta) \afun(\theta).
\]
We emphasize that the cohomological equation for $\afun$, 
\begin{equation}
\afun(\theta) - \sigma\afun(\theta+\omega)
= -t_0(\theta). 
\label{eq:Cohom_afun_NT}
\end{equation}
can be solved in Fourier space: 
\[
 \afun(\theta)= \sum_{k\in\ZZ} \frac{-t_{0k}}{1-\sigma e^{2\pi{\ii}k\omega}}  e^{2\pi{\ii}k\theta}.
\]
(Notice that, since $|\sigma|<1$, the divisors are uniformly far from $0$.)

Second, we write the correction term of the parameterization of the rotational circle as 
$\Delta K(\theta)=  P(\theta)\xi(\theta)$,
where $\xi:\TT\to \RR^2$ is a periodic function.  
Then, by multiplying \eqref{eq:Newton_qp}  by $P(\theta+\omega)^{-1}$,  using approximate reducibility \eqref{approximate_reducibility_qp} 
and neglecting quadratically small terms, we obtain the following cohomological equation 
\begin{equation}
\label{eq:Newton_qp}
\Lambda(\theta)\xi(\theta)-\xi(\theta+\omega) + B_a(\theta)\Delta a + B_\mu(\theta)\Delta\mu 
= \eta(\theta),
\end{equation}
where 
$B_a(\theta)= P(\theta+\omega)^{-1} \Dif_a F(K(\theta),a,\mu,\eps_0)$,  
$B_\mu(\theta)= P(\theta+\omega)^{-1} \Dif_\mu F(K(\theta),a,\mu,\eps_0)$, 
and $\eta(\theta)= -P(\theta+\omega)^{-1} E(\theta)$ is the error of invariance  in the adapted frame. 

Notice that the previous system is  diagonal, and it splits into
\begin{eqnarray}
\label{eq:cohom_tan_qp}
\xiL(\theta) -\xi^\tL(\theta+\omega) + B^\tL_a(\theta)\Delta a + B^\tL_\mu(\theta)\Delta\mu &=& \etaL(\theta),\\
\label{eq:cohom_stab_qp}
\sigma\xiN(\theta) -\xi^\tN(\theta+\omega) + B^\tN_a(\theta)\Delta a + B^\tN_\mu(\theta)\Delta\mu &=& \etaN(\theta),
\end{eqnarray}
where
\[
\begin{split}
	B^\tL_a(\theta) = \phantom{-}N(\theta+\omega)^\top \Omega \ \Dif_a F(K(\theta),a,\mu), &\quad
	B^\tL_\mu(\theta) =\phantom{-} N(\theta+\omega)^\top \Omega \ \Dif_\mu F(K(\theta),a,\mu) \\
	B^\tN_a(\theta) = -L(\theta+\omega)^\top \Omega \ \Dif_a F(K(\theta),a,\mu), & \quad
	B^\tN_\mu(\theta) = -L(\theta+\omega)^\top \Omega \ \Dif_\mu F(K(\theta),a,\mu).
\end{split}
\]
In particular: $b_a(K,a,\mu,\eps_0)= \avg{B^\tL_a(\theta)}$, $b_\mu(K,a,\mu,\eps_0)= \avg{B^\tL_\mu(\theta)}$.

It is the moment to face cohomological equations \eqref{eq:cohom_tan_qp} and \eqref{eq:cohom_stab_qp}, 
which are in fact very different, and introduce some notation.  We will denote by $\xi= \calR_\sigma \eta$ the solution of 
\[
	 \sigma\xi(\theta) -\xi(\theta+\omega)  = \eta(\theta),
\]
that is, in Fourier series:
\[
 \xi(\theta)= \calR_\sigma\eta(\theta)= \sum_{k\in\ZZ} \frac{\eta_{k}}{\sigma-e^{2\pi{\ii}k\omega}}  e^{2\pi{\ii}k\theta}.
\]
Notice again that, since $|\sigma|<1$, the divisors are uniformly far from $0$. The case $\sigma=1$ is very different: 
the right hand side has to have zero average, the solution if exists it is not unique, and the divisors can be arbitrarily small.    
We will denote by $\xi= \calR\eta$ the solution of 
\[
	 \xi(\theta) -\xi(\theta+\omega)  = \eta(\theta) - \avg{\eta}
\]
with zero average, 
that is, in Fourier series:
\[
 \xi(\theta)= \calR\eta(\theta)= \sum_{k\in\ZZ^*} \frac{\eta_{k}}{1-e^{2\pi{\ii}k\omega}}  e^{2\pi{\ii}k\theta}.
\]
The solution involves small divisors and it suffices Diophantine conditions on $\omega$ to ensure the convergence 
of the expansions. Notice also the adjustment to get zero average in the right hand side of the small divisors equation, 
and that we can add constants to $\xi$ to get (non-zero average) solutions. 

Third, for $\Delta a$ fixed, we solve  \eqref{eq:cohom_tan_qp} and \eqref{eq:cohom_stab_qp} as follows.
We compute $\Delta\mu= \Delta\mu[\Delta a]$ by adjusting averages in \eqref{eq:cohom_tan_qp}, so that 
\[
	\Delta\mu= \frac{\avg{\etaL} - \avg{B^\tL_a}\Delta a}{\avg{B^\tL_\mu}} \simeq
	\frac{\avg{\etaL} - b_a^0 \Delta a}{b_\mu},
\]
where in the last approximation we are skipping second order error terms. Notice that we need a twist condition with 
respect to the adjusting parameter $\mu$. We emphasize the dependence of $\Delta\mu$ on $\Delta a$ (as we will do in the sequel for other objects). 
With this choice of $\Delta\mu$ we compute $\xiL= \xiL[\Delta a]$, $\xiN= \xiN[\Delta a]$ as follows:
\[
\xiN(\theta)= \calR_\sigma\eta(\theta) - \calR_\sigma B^\tN_a(\theta)\ \Delta a- \calR_\sigma B^\tN_\mu(\theta)\ \Delta \mu
\]
for the solution of \eqref{eq:cohom_stab_qp}
\[
\hat\xiL(\theta)= \calR\eta(\theta) - \calR B^\tN_a(\theta)\ \Delta a- \calR B^\tN_\mu(\theta)\ \Delta \mu
\]
for the zero-average solution of \eqref{eq:cohom_tan_qp}, 
\[
\xiL_0= -e_p-\avg{L^x( \theta) \hat\xiL(\theta) + N^x(\theta) \xiN(\theta)}
\]
to fix the phase (notice that $\avg{L^x}= 1$) and, finally 
\[
\xiL(\theta)= \xiL_0 + \hat\xiL(\theta).
\] 

Fourth, we obtain a correction $\Delta K= \Delta K[\Delta a]$ for improving \eqref{eq:InvarianceError} and 
\eqref{eq:PhaseError}:
\[
 \Delta K[\Delta a](\theta)= L(\theta) \xiL[\Delta a](\theta) + N(\theta) \xiN[\Delta a](\theta).
\]

In summary, from the previous four steps we obtain a univariate function 
\[
	\Delta a \to b_a[\Delta a]= b_a(K+\Delta K[\Delta a], a+\Delta a, \mu +\Delta \mu[\Delta a], \eps_0)
\]
for which we have to solve the equation 
\begin{equation}
\label{eq:Delta_a equation}
	b_a[\Delta a] - b_a^0 = 0.
\end{equation}
In the implementation of each step of Newton method, instead of solving this equation, we apply one step of
Steffensen's method to this equation starting with $\Delta a= 0$.  In the implementation, we control the
non-degeneracy condition to 
solve \eqref{eq:Delta_a equation}. 

With the previous Newton method we compute an invariant rotational circle with fixed $a$-twist
for a fixed value of $\eps_0$. 
In order to implement the continuation with respect to parameter $\eps$ one can compute derivatives of
$(K,a,\mu)$ with 
respect to $\eps$, at $\eps_0$. The type of equations one has to solve  are of the
same type as to perform a Newton step.
In particular, one has  \eqref{eq:Newton_qp} with
 \[
E(\theta)= \frac{\partial F}{\partial \eps}(K(\theta),a,\mu,\eps_0),
\]
and $\Delta K= \frac{\partial K}{\partial\eps}$, $\Delta a= \frac{\partial a}{\partial\eps}$, 
$\Delta_\mu= \frac{\partial \mu}{\partial \eps}$.

For the implementation of Newton's method and continuation method described here we use
Fourier series to 
represent periodic functions. Thus, we use FFTs  to switch from grid representation to Fourier representation. 
All operations can be done at linear cost in grid or Fourier representations, except the ones switching
representations. 
Hence, the cost of the algorithms is $O(N \log(N))$ where $N$ is the size of the representation
(the size of the grid or the number of Fourier modes). See e.g. \cite{CallejaL09,
CallejaL10, HaroCFLM16}
for some guidelines.

\section{The applications}
\label{section:applications}

In this section, we implement the algorithms presented in this paper for some specific families of conformally symplectic maps with conformal factor 
$\sigma$, referred to as 
 \emph{dissipative standard non-twist maps}. These are 
defined by the dynamical systems $F_{a,\mu,\eps}:\TT\times\RR \to \TT\times\RR$ given by
\begin{equation}
F_{a,\mu,\eps}\begin{pmatrix}x\\ y\end{pmatrix}
=
\begin{pmatrix}
x + \left(\sigma y + {\eps} p(x)- a \right)^2+\mu\\
\sigma y +  {\eps} p(x)
\end{pmatrix},
\label{eq:DSNTM}
\end{equation}
where $p:\TT\to \RR$ is a $1$-periodic function, $a$, $\mu$ are adjusting parameters (whose roles will be explained below), and $\eps$ is the perturbative parameter. 
We will consider two examples: (symmetric) $p(x)= \frac{1}{2\pi} \sin(2\pi x)$, and (non-symmetric) 
$p(x)= \frac{1}{2\pi} \left(\sin(2\pi x) + \cos(4\pi x)\right)$. The 
names we give to the functions $p(x)$ will be justified later.

\subsection{Preliminaries}

We start by analyzing \eqref{eq:DSNTM} for $\eps=0$, which is integrable. For each $a,\mu$, 
there is an invariant circle parameterized by 
\begin{equation*}
K_{a,\mu,0}(\theta)=
\begin{pmatrix}
\theta\\
0
\end{pmatrix}, 
\end{equation*}
whose internal  dynamics is given explicitly by 
\begin{equation*}
f_{a,\mu,0}(\theta)=
\theta + a^2 + \mu.
\end{equation*}
Moreover, the adapted frame and the corresponding linearized dynamics are
\begin{equation*}
P(\theta)
=\begin{pmatrix}\Dif K(\theta) & N(\theta)\end{pmatrix} 
= \begin{pmatrix}1&0\\ 0&1\end{pmatrix},
\hspace{0.5cm}
\Lambda(\theta)
= \begin{pmatrix}\Dif f(\theta)&0\\ 0&\LambdaN(\theta)\end{pmatrix}
= \begin{pmatrix}1&0\\ 0& \sigma\end{pmatrix}.
\end{equation*}

If we are looking for an initial invariant tori with fixed quasi-periodic frequency $\omega$, then parameter $a,\mu$ are linked by the relation   
\begin{equation}
\label{eq:firstfreq}
\omega = a^2 + \mu.
\end{equation}
That is, $\mu= \omega - a^2$.

The $a$-twist is 
\[
    b_a(K; a,\mu,0)= \avg{N(\theta+\omega)^\top \Omega\ \Dif_a F_{a,\mu,0}(K(\theta))}= 2a,
\]
while the $\mu$-twist is 
\[
     b_\mu(K; a,\mu,0)= \avg{N(\theta+\omega)^\top \Omega\ \Dif_\mu F_{a,\mu,0}(K(\theta),a,\mu)}= 1.
\]
Since the  $\mu$-twist is non zero  we can isolate $\mu$.
Notice however, that the invariant circle 
with frequency $\omega$ is non-$a$-twist for whenever $a= 0$ (and hence $\mu= \omega$).

Since the $\mu$-twist is non zero, from an  implicit function theorem we get that 
for $a,\eps$ close to zero, we can find $\mu= \mu(a,\eps)$ and a circle parameterized by $K_{a,\mu,\eps}$ which 
is invariant for $F_{a,\mu,\eps}$ and whose internal dynamics is a rotation with frequency $\omega$. 
So given $a,\eps$, the unfolding parameter $\mu$ is used to adjust the frequency to $\omega$.
By writing $\bar b_a(a,\eps)= b_a(K_{a,\mu(a,\eps),\eps}; a,\mu(a,\eps),\eps)$, then the equation
we need to solve is
\[
	\bar b_a(a,\eps)= 0, 
\]
and to apply the implicit function theorem in order to find $a$ for small
enough $\eps$ we also need that 
\[
	\frac{\partial \bar b_a}{\partial a} (0,0) \neq 0.
\]
In our example, $\frac{\partial \bar b_a}{\partial a} (0,0)= 2.$

\subsection{Continuation of  the non-$a$-twist circle in the symmetric case}

In this section we study the family \eqref{eq:DSNTM} with $p(x)= \tfrac{1}{2\pi} \sin(2\pi x)$. Since 
$p(x-\tfrac12)= -p(x)$ then the involution $S:\TT\times \RR \to \TT\times \RR$ defined by 
\[
	S\begin{pmatrix} x \\ y \end{pmatrix} = \begin{pmatrix} x-\tfrac12 \\ -y \end{pmatrix}
\]
is a symmetry of the family with respect to parameter $a$, meaning that 
\begin{equation*}
	S \circ F_{a,\mu,\eps} \circ S =  F_{-a,\mu,\eps}.
\end{equation*}
This symmetry property implies that if $K_{a,\mu,\eps}$ is a parameterization of an invariant circle
for $F_{a,\mu,\eps}$ with internal dynamics $f_{a,\mu,\eps}$, then $K_{-a,\mu,\eps}= S \circ K_{a,\mu,\eps}$
is  parameterization of 
an invariant circle for $F_{-a,\mu,\eps}$ with internal dynamics $f_{-a,\mu,\eps}= f_{a,\mu,\eps}$.
In particular, for $a= 0$, the invariant circle parameterized by $K_{0,\mu,\eps}$ is $S$-symmetric, since it is also 
parameterized by $S \circ K_{0,\mu,\eps}$. In fact, since they is a unique parameterization such that 
$\avg{K^x_{0,\mu,\eps}(\theta)-\theta}= 0$, it satisfies:
\[
	K_{0,\mu,\eps}(\theta)= S \circ K_{0,\mu,\eps}(\theta+\tfrac12).
\]
In this case, by selecting $\mu$ so that $f_{0,\mu,\eps}(\theta)= \theta+\omega$, we have
\[
	\bar b_a(0,\eps)= b_a(K_{0,\mu,\eps}, 0, \mu, \eps)= 0.	
\]
That is, the tori with $a= 0$ are non-$a$-twist. This example will be a first test of our algorithms. 

Below we will explain the results derived from our implementation of algorithms in sections \ref{sec:NT-tori}
and \ref{sec:GeneralConf} to this example, for the dissipative parameter $\sigma=0.8$. 

\begin{figure}[ht]
\begin{tabular}{cc}
\includegraphics[width=0.35\linewidth]{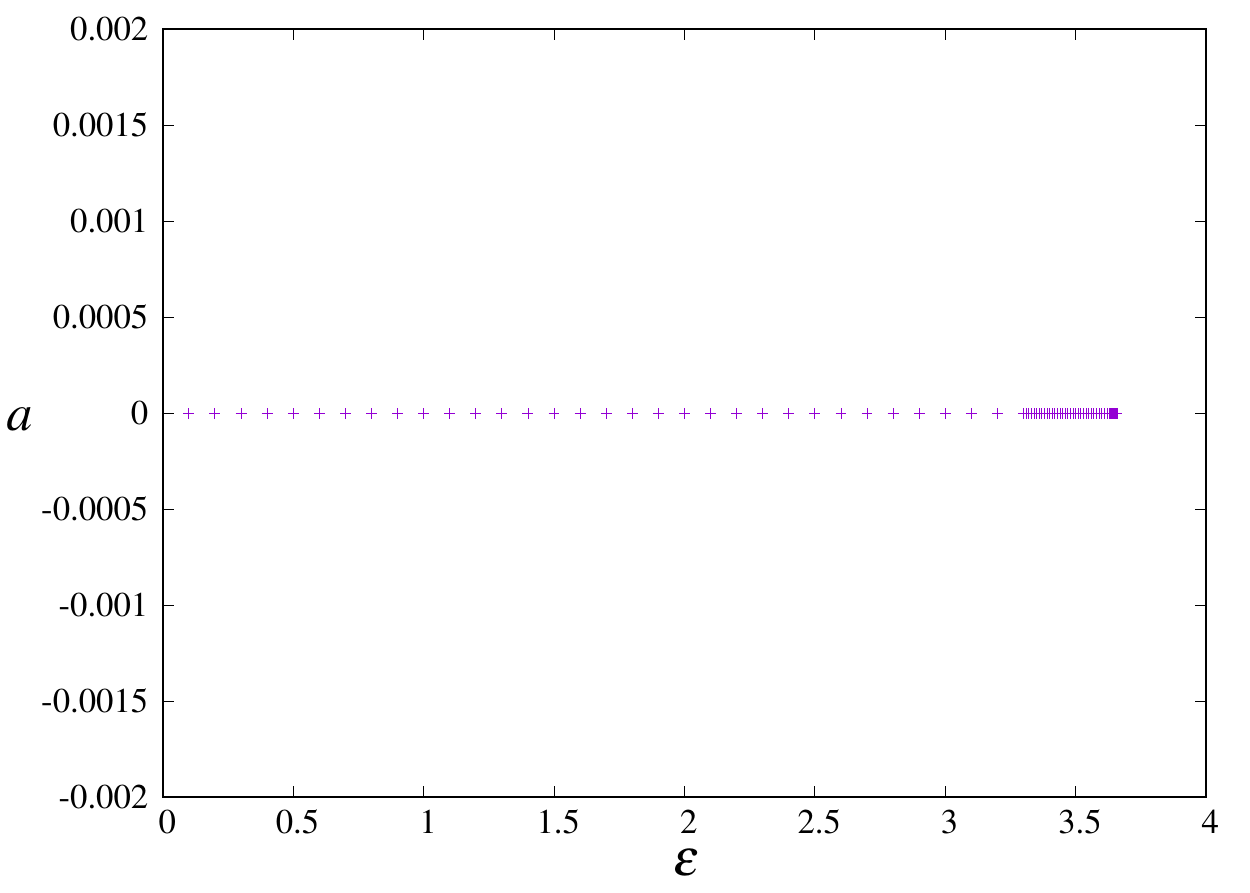} &
\includegraphics[width=0.35\linewidth]{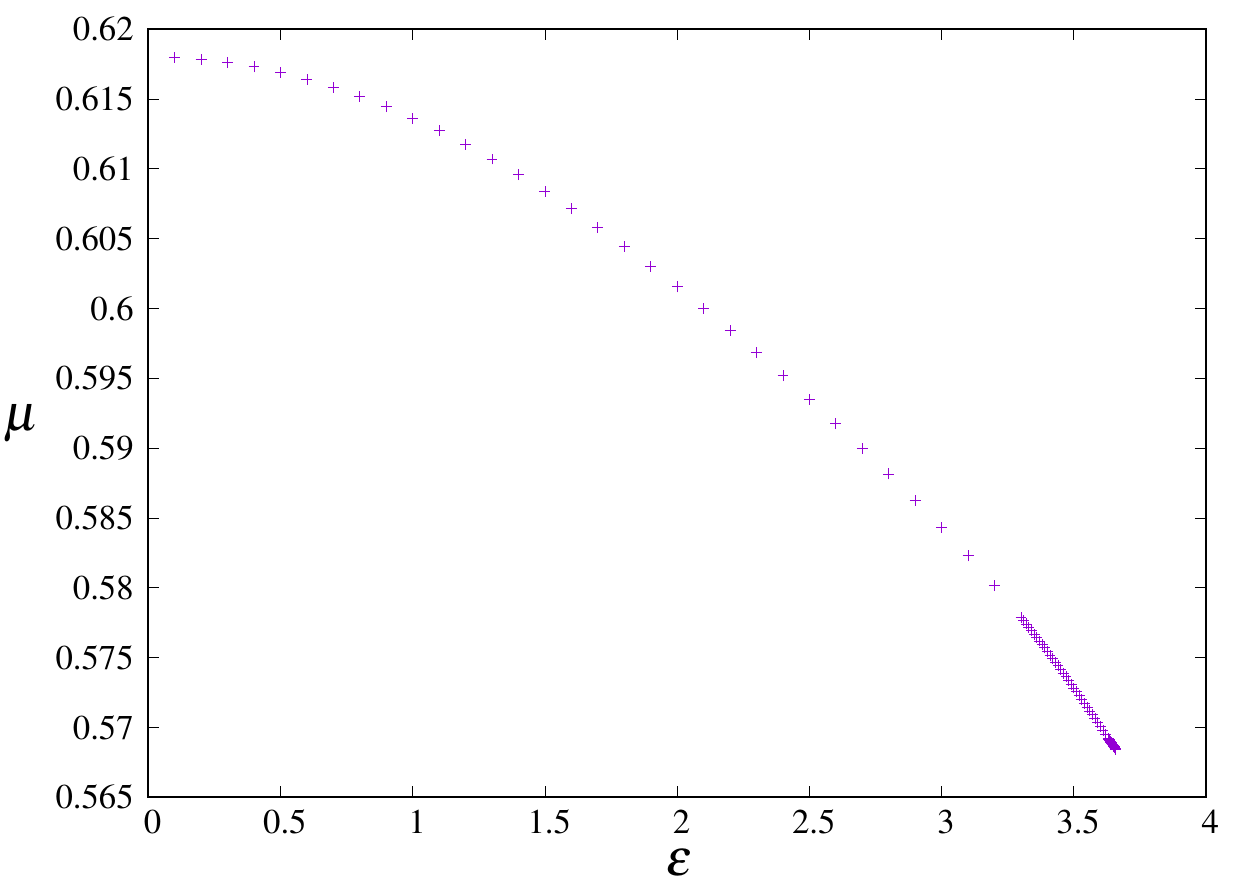}
\end{tabular}
\caption{\label{fig: ex1 parameters} 
Continuation w.r.t. $\eps$ of a non-$a$-twist circle with frequency $\omega$ (symmetric case): 
(left) adjusting parameter $a$; (right) unfolding parameter $\mu$.
}
\end{figure}

\begin{figure}
\begin{subfigure}[a][$\eps= 2.000000$, $a=  0.000000$, $\mu= 0.6015602$]{ 
\centering
\begin{tabular}{cc}
\includegraphics[width=0.35\linewidth]{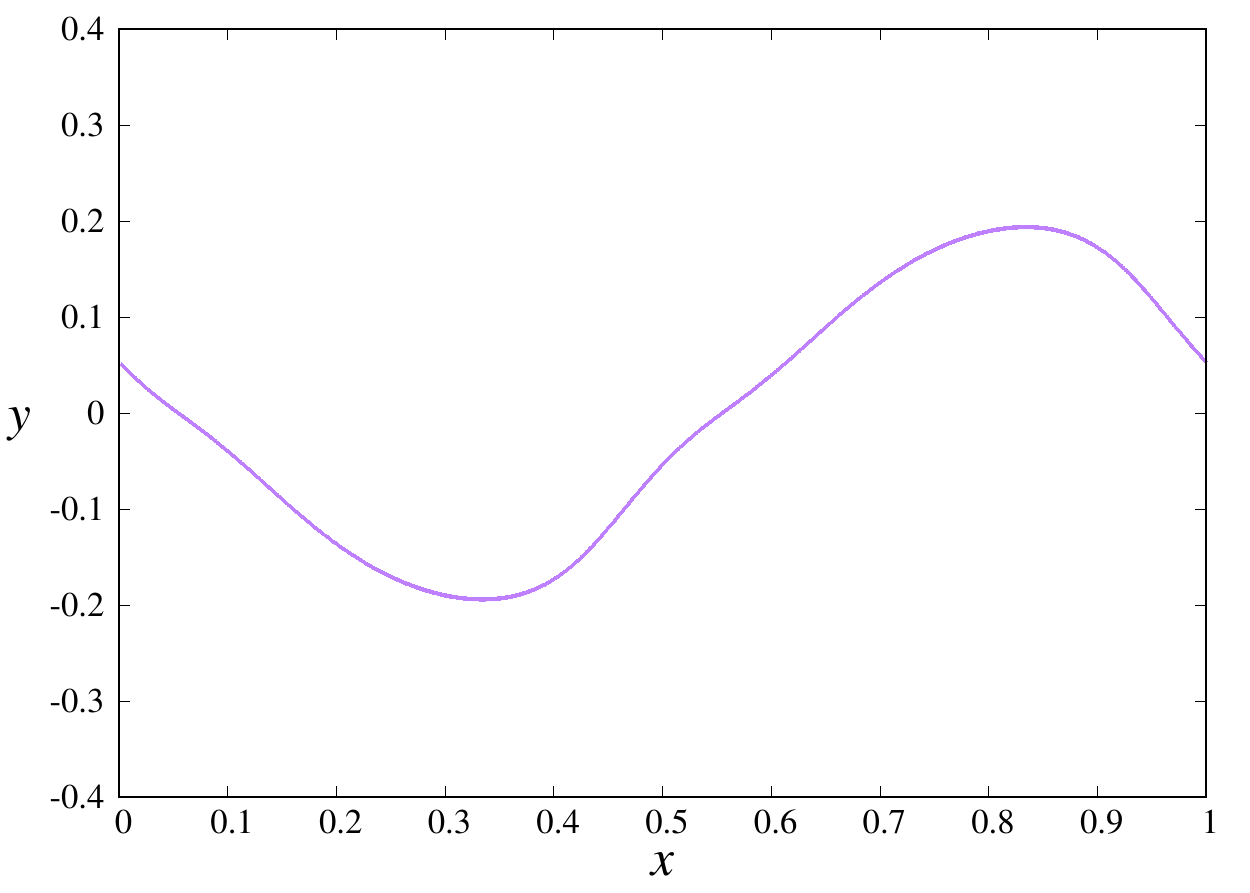} 
&
\includegraphics[width=0.35\linewidth]{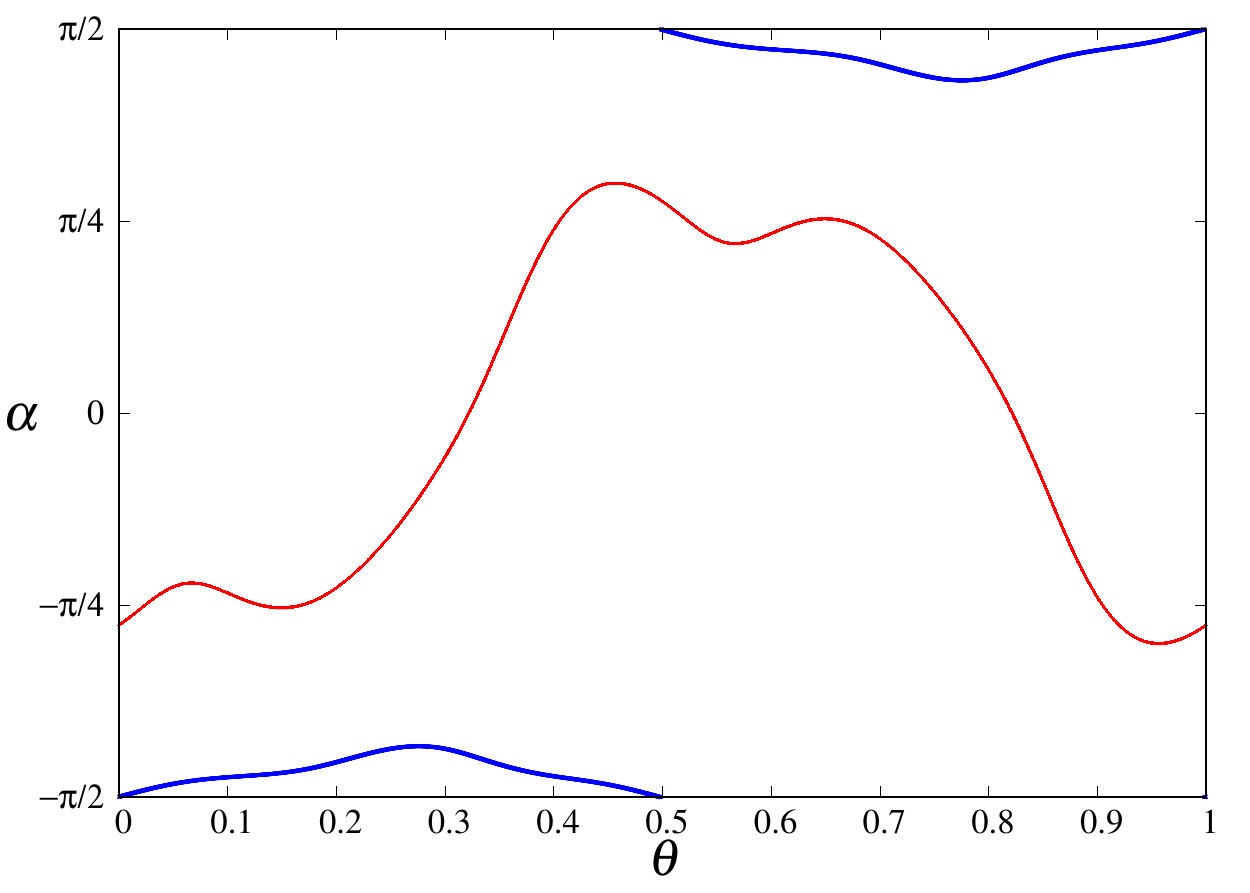} 
\end{tabular}
}
\end{subfigure}
\begin{subfigure}[a][$\eps= 3.000000$, $a= 0.000000$, $\mu= 0.5843217$]{
\centering
\begin{tabular}{cc}
\includegraphics[width=0.35\linewidth]{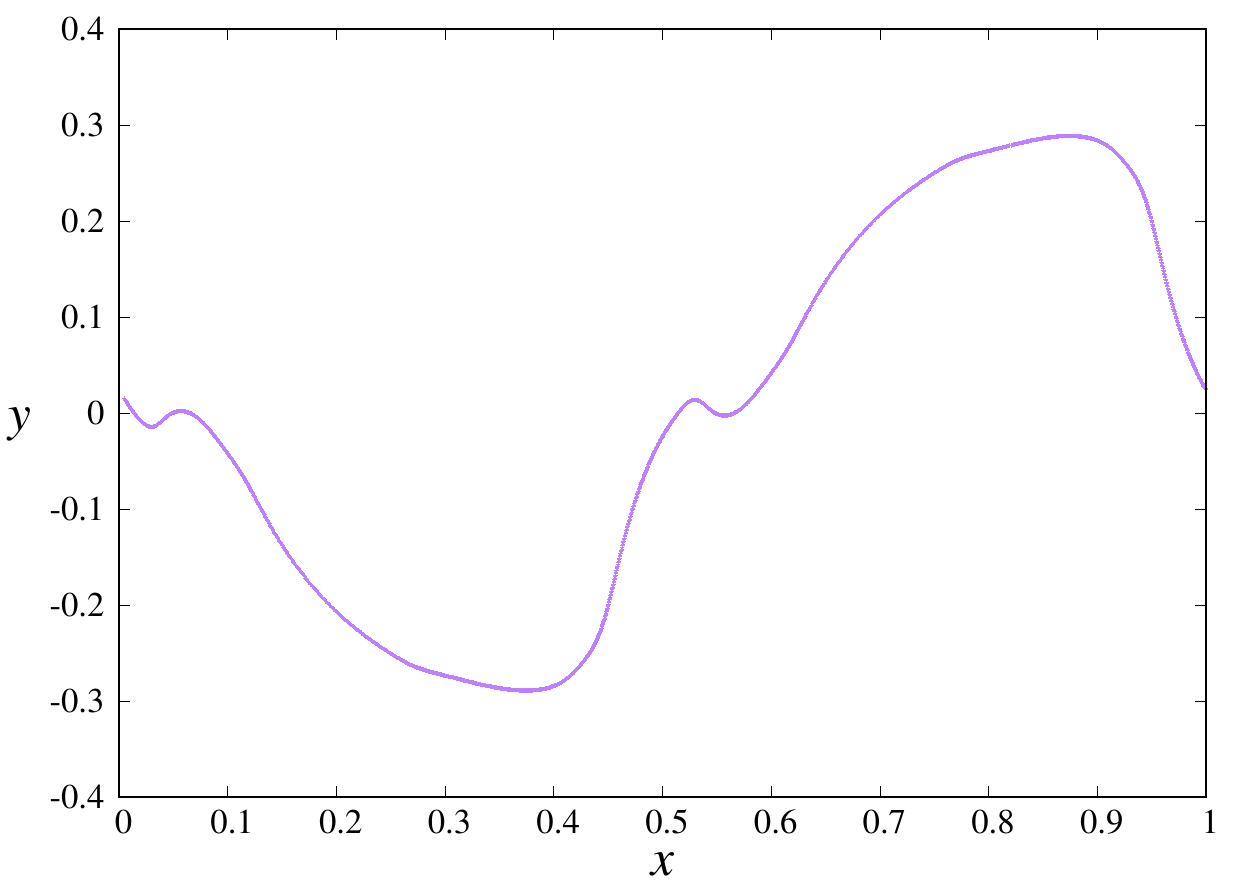} 
&
\includegraphics[width=0.35\linewidth]{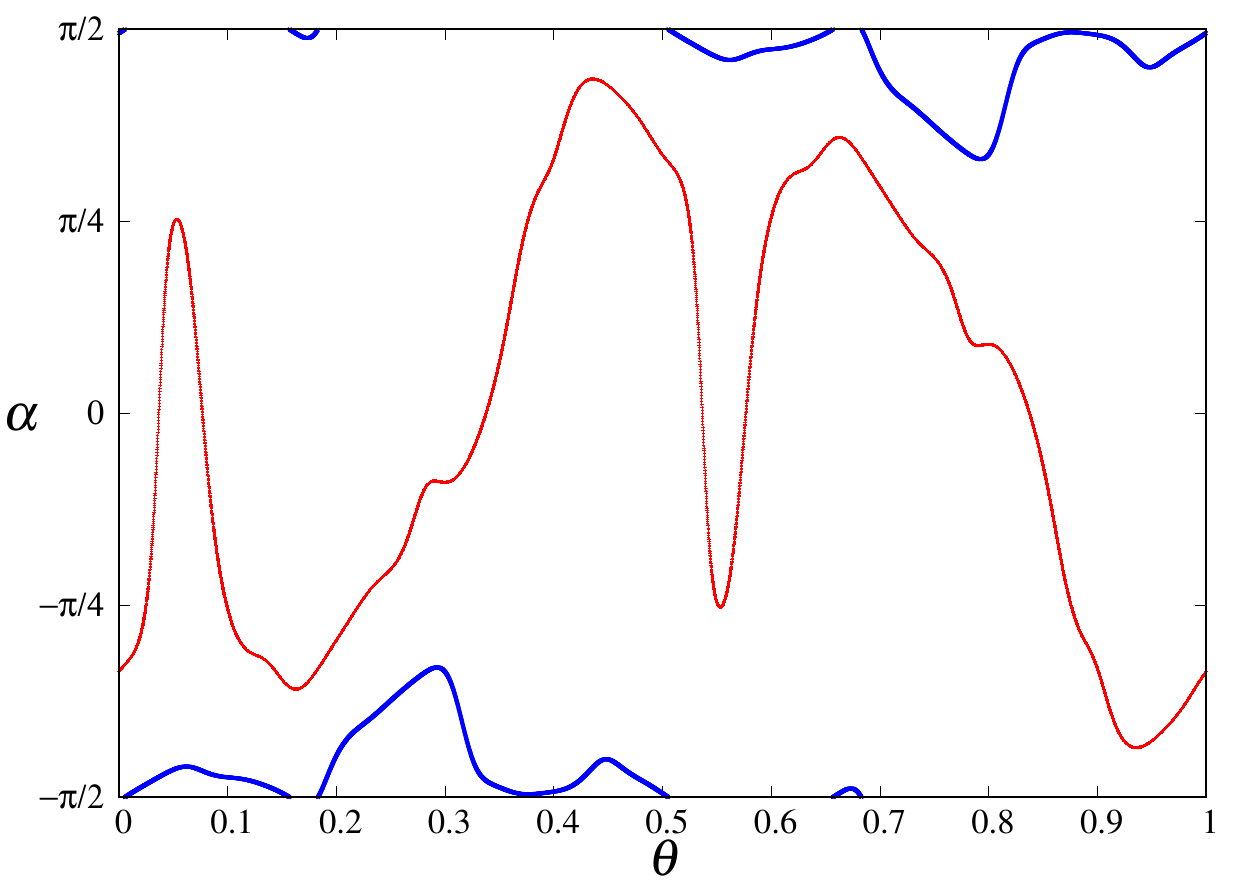} 
\end{tabular}
}
\end{subfigure}
\begin{subfigure}[a][$\eps= 3.658600$, $a= 0.000000$, $\mu= 0.5684363$]{
\centering
\begin{tabular}{cc}
\includegraphics[width=0.35\linewidth]{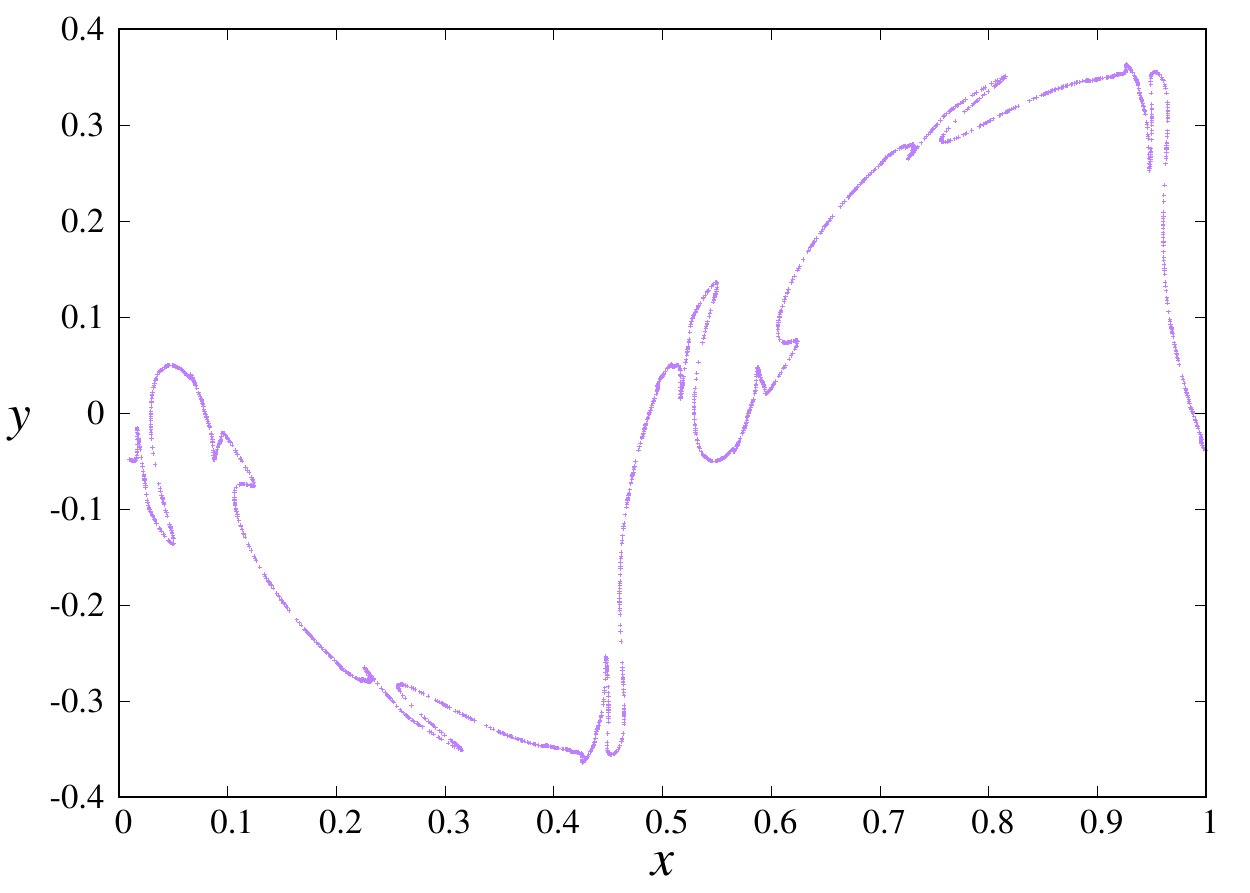} 
& 
\includegraphics[width=0.35\linewidth]{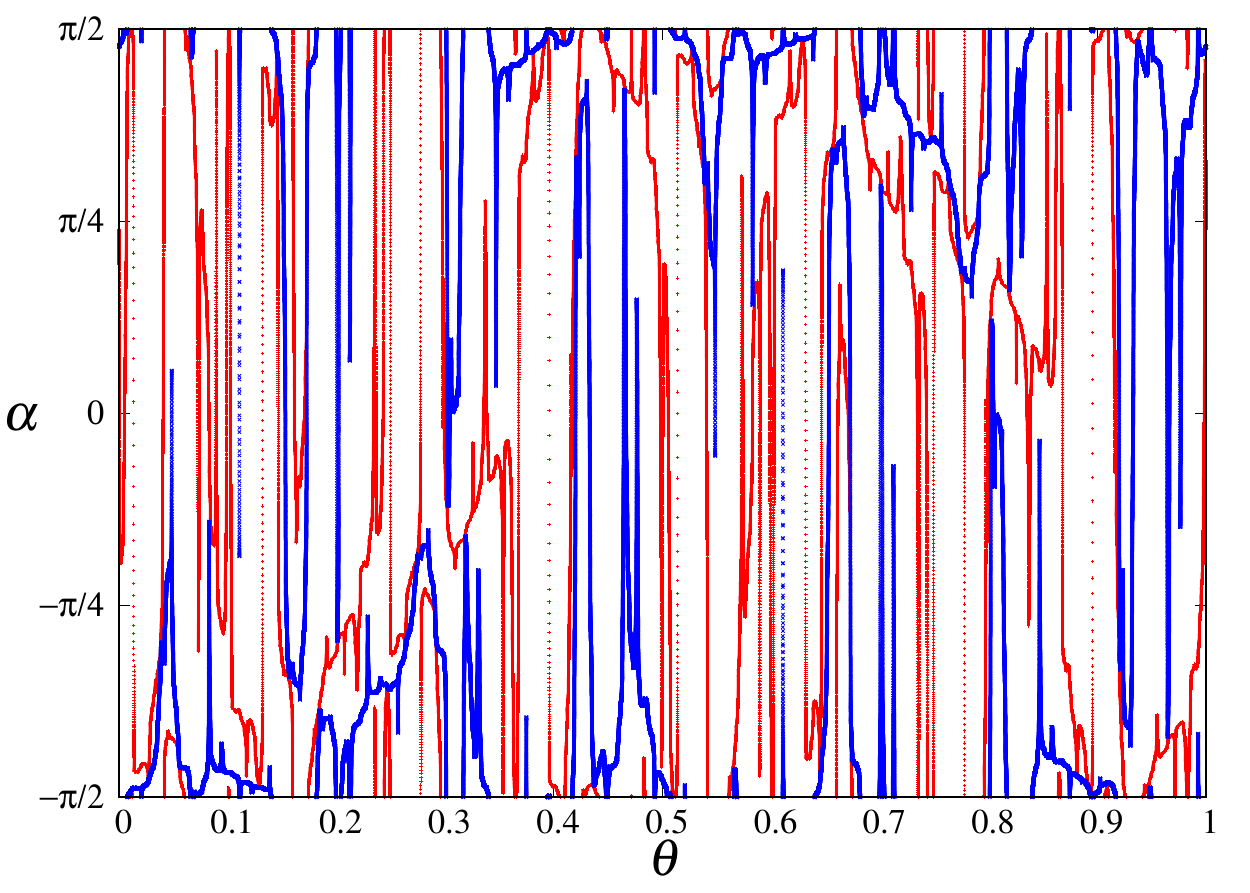}
\end{tabular}
}
\end{subfigure}
\caption{\label{fig: ex1 non-twist circle}
Continuation w.r.t. $\eps$ of a non-$a$-twist circle with frequency $\omega$ (symmetric case):
(left) invariant circle; (right) projectivized tangent bundle (in red) and stable bundle (in blue).
}
\end{figure}

\begin{figure}[ht]
\begin{tabular}{cc}
\includegraphics[width=0.35\linewidth]{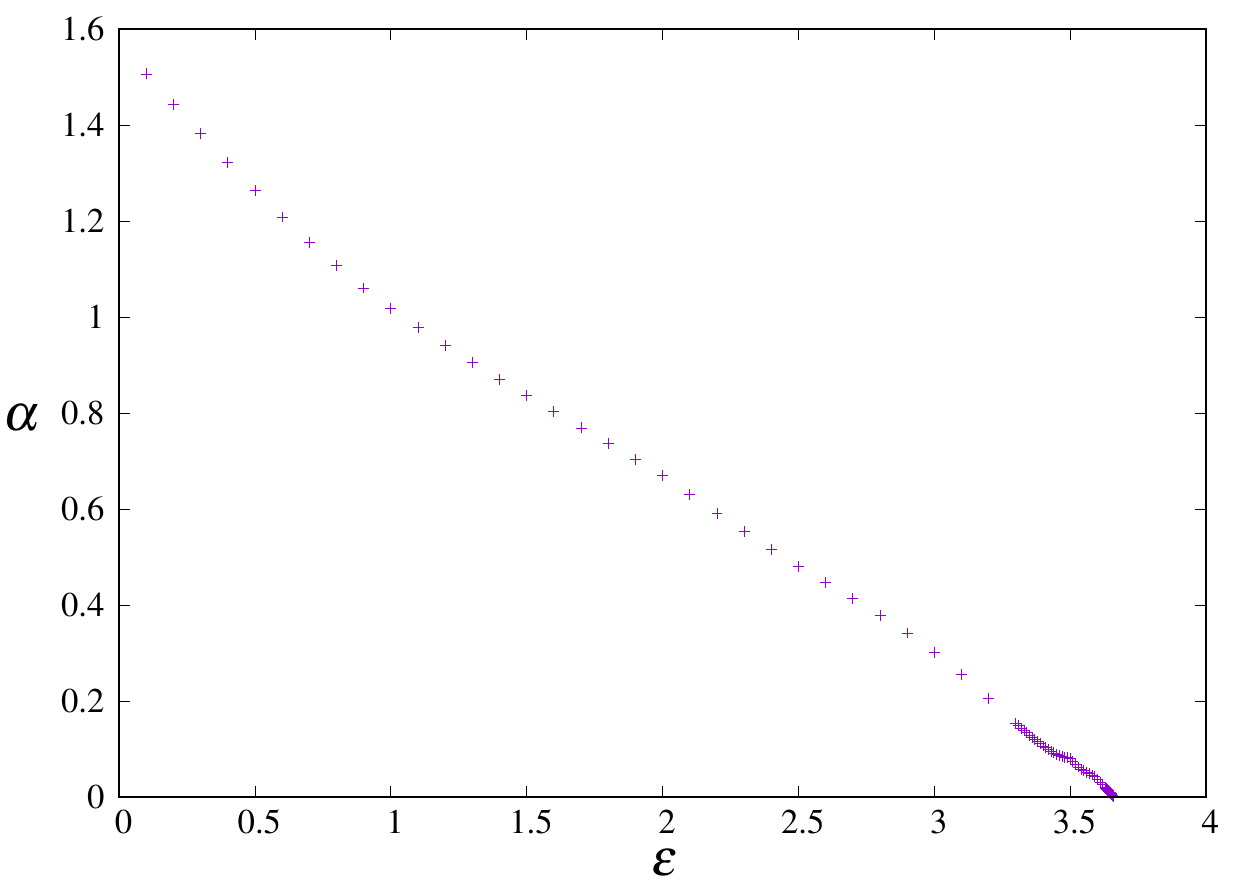} &
\includegraphics[width=0.35\linewidth]{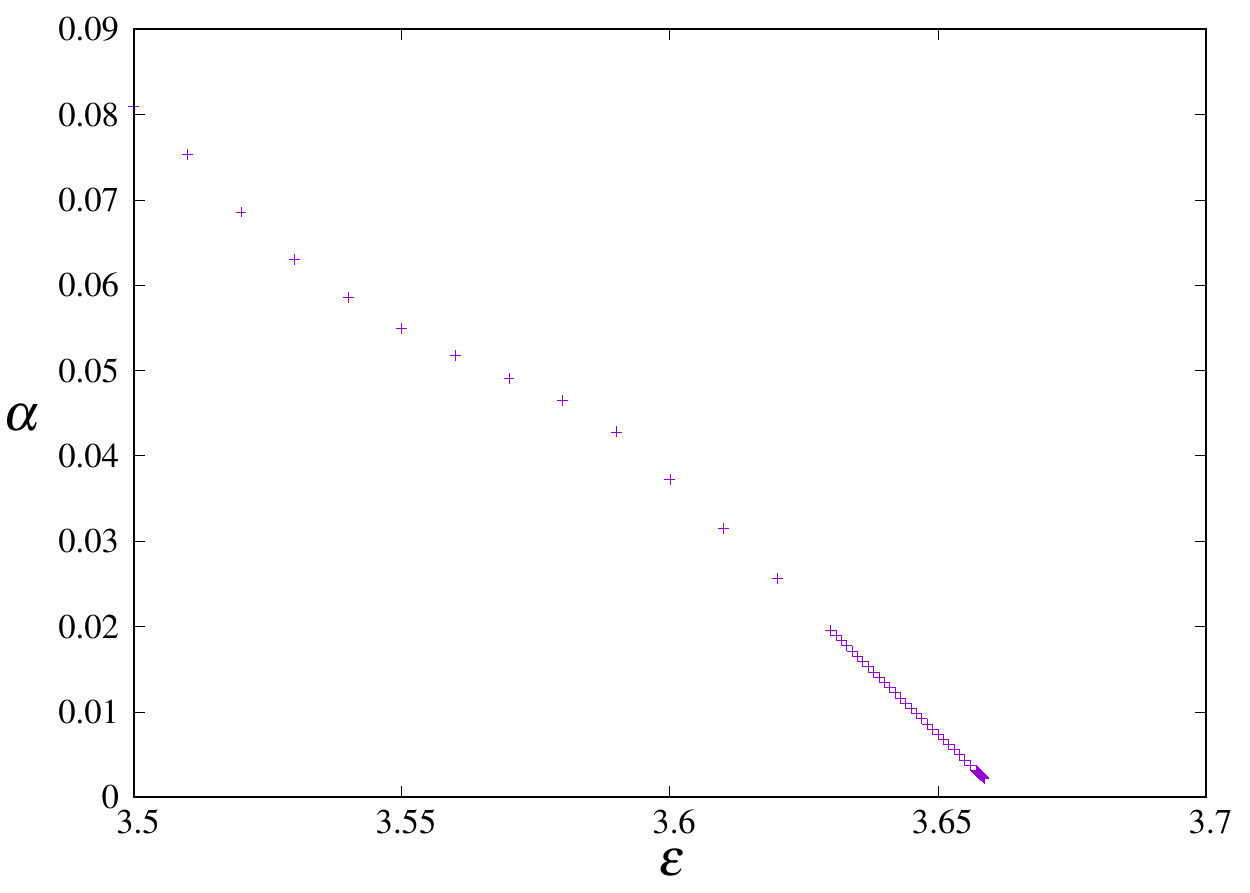}
\end{tabular}
\caption{\label{fig: ex1 angles}
Continuation w.r.t. $\eps$ of a non-$a$-twist circle with frequency $\omega$ (symmetric case):
(left) minimum angle $\alpha$ between the stable and tangent bundles as function of  $\eps$;
(right) critical behavior. The breakdown of the circle is produced at $\eps_{\rm c}\simeq 3.662396$.              
}
\end{figure}

Using algorithm of Section~\ref{sec:NT-tori}, we continue with respect to parameter $\eps$ a non-$a$-twist circle 
with rotation number $\omega= \tfrac{1}{2} (\sqrt{5}-1)$, and adjust parameters $a$ and $\mu$ accordingly. 
The adjusting parameters are shown in Figure~\ref{fig: ex1 parameters}.
Starting from $\eps= 0$, the continuation goes up to $\eps= 3.658600$, close to breakdown, in which the number 
of Fourier modes demanded by the algorithm is $262144$. Some of the non-$a$-twist circles are shown in 
Figure~\ref{fig: ex1 non-twist circle}, together with the corresponding tangent and stable bundles, represented 
by their angles with respect to the horizontal axis $\alpha$. The complex behavior observed in the bundles preludes 
the breakdown of the invariant circle. Notice that when both bundles collide, the normal hyperbolicity property fails, 
and this happens even though the contraction factor is far from $1$ (it is $\sigma= 0.8$).
This collision behavior has been observed in other contexts 
\cite{CanadellH14,CanadellH17b,FiguerasH15,HaroL06c,HaroL07}, and in \cite{CallejaF12} for 
$a$-twist circles in conformally symplectic systems. From these references one conjectures
that, even though the behavior is very wild, there is some sort of regularity and the minimum 
angle between the invariant bundles behaves very smoothly, 
in fact asymptotically in a linear fashion when approaching the breakdown, as shown in 
Figure~\ref{fig: ex1 angles}. This behavior lets us extrapolate the critical breakdown
parameter very consistingly,
being $\eps_{\rm c}\simeq 3.662396$. 

The symmetry properties of the family lead to 
several features. First, parameter $a$ is always $0$, as it is shown 
in Figure~\ref{fig: ex1 parameters}. Moreover, the non-$a$-twist circles and their bundles have also 
symmetry properties, as it is shown in Figure~\ref{fig: ex1 non-twist circle}. In particular, 
We note that the collapse in this symmetric case happens on both sides of the bundles.
We expect that when the bundles collapse, there will be no gap between
the bundles on either side of the bundles with respect to $\alpha$. Later in
this section, we will see that in the nonsymmetric 
the bundles collapse leaving a gap between the bundles for all values
of $\theta$, but only on one side of the bundles with respect to $\alpha$.

\begin{figure}
\begin{tabular}{cc}
\includegraphics[width=0.35\linewidth]{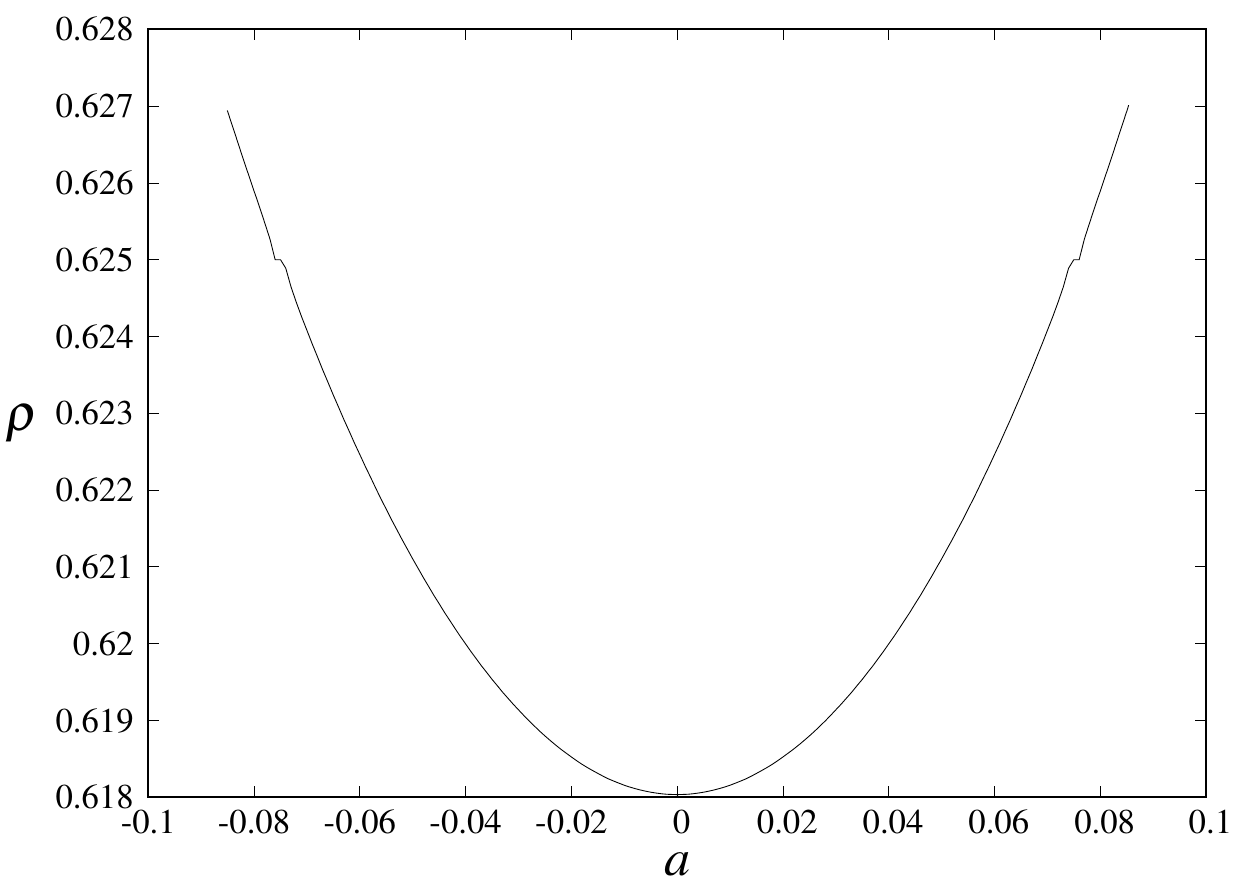} &
\includegraphics[width=0.35\linewidth]{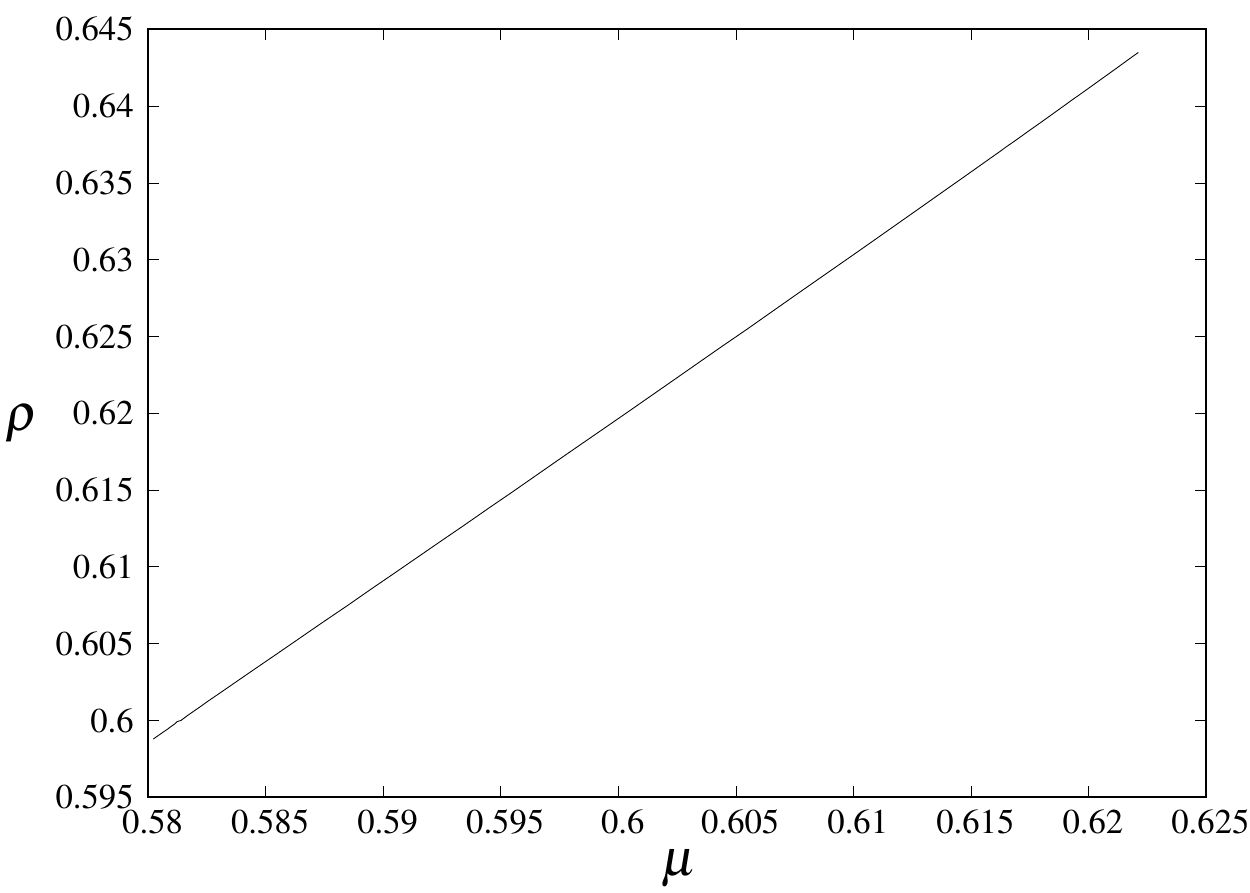}
\end{tabular}
\caption{\label{fig:ex 1 rotation number}
Rotation number versus continuation parameter from the non-$a$-twist circle 
with $a= 0.000000$, $\mu= 0.5984626$, $\eps= 2.20000$ (symmetric case):\newline
(left) continuation w.r.t. $a$; (right) continuation w.r.t. $\mu$.}
\end{figure}

We also performed some computations to illustrate that the analytic condition that the invariant
circle is non-$a$-twist translates into dynamical properties of the rotation number of the invariant
cicle when we move parameters.
We implement the algorithm in Section~\ref{sec:GeneralConf} to continue invariant tori regardless the internal
dynamics and compute the corresponding rotation number,
by starting with a non-$a$-twist circle from the previous implementation. In particular, 
we have selected a non-$a$-twist circle for $\eps= 2.2$, so that $\mu=  1.5984626393$ and $a= 0$. 
We first perform continuations for $\mu$ and $\eps$ fixed, increasing and decreasing the parameter $a$, respectively. 
The graph of the rotation number of the invariant circle as a function of $a$ is shown in Figure~\ref{fig:ex 1 rotation number} (Left). As expected, 
the non-$a$-twist circle corresponds to a critical point of this graph. The graph is symmetric, also as expected 
from the symmetry properties of the family being studied. 
Notice also the presence of visible resonances, corresponding to rotation number $5/8$. 
However, by performing a continuation with respect to $\mu$ instead of $a$
(and starting with the same initial torus), we observe that the starting
torus does not correspond to a minimum of the rotation number as a function of $\mu$, as shown in 
Figure~\ref{fig:ex 1 rotation number} (Right).
This is because the non-twist-property is associated to parameter $a$, and the invariant circle is $\mu$-twist. 

\subsection{Continuation of  the non-$a$-twist circle in the nonsymmetric case}

In this section we consider the family \eqref{eq:DSNTM} with $p(x)= \tfrac{1}{2\pi} (\sin(2\pi x)+\cos(4\pi x)$, that 
(apparently) does not have symmetry properties. We again
take $\sigma= 0.8$, and $\omega= \tfrac12(\sqrt{5}-1)$.
We have followed the same plan as in previous example.

\begin{figure}[ht]
\begin{tabular}{cc}
\includegraphics[width=0.35\linewidth]{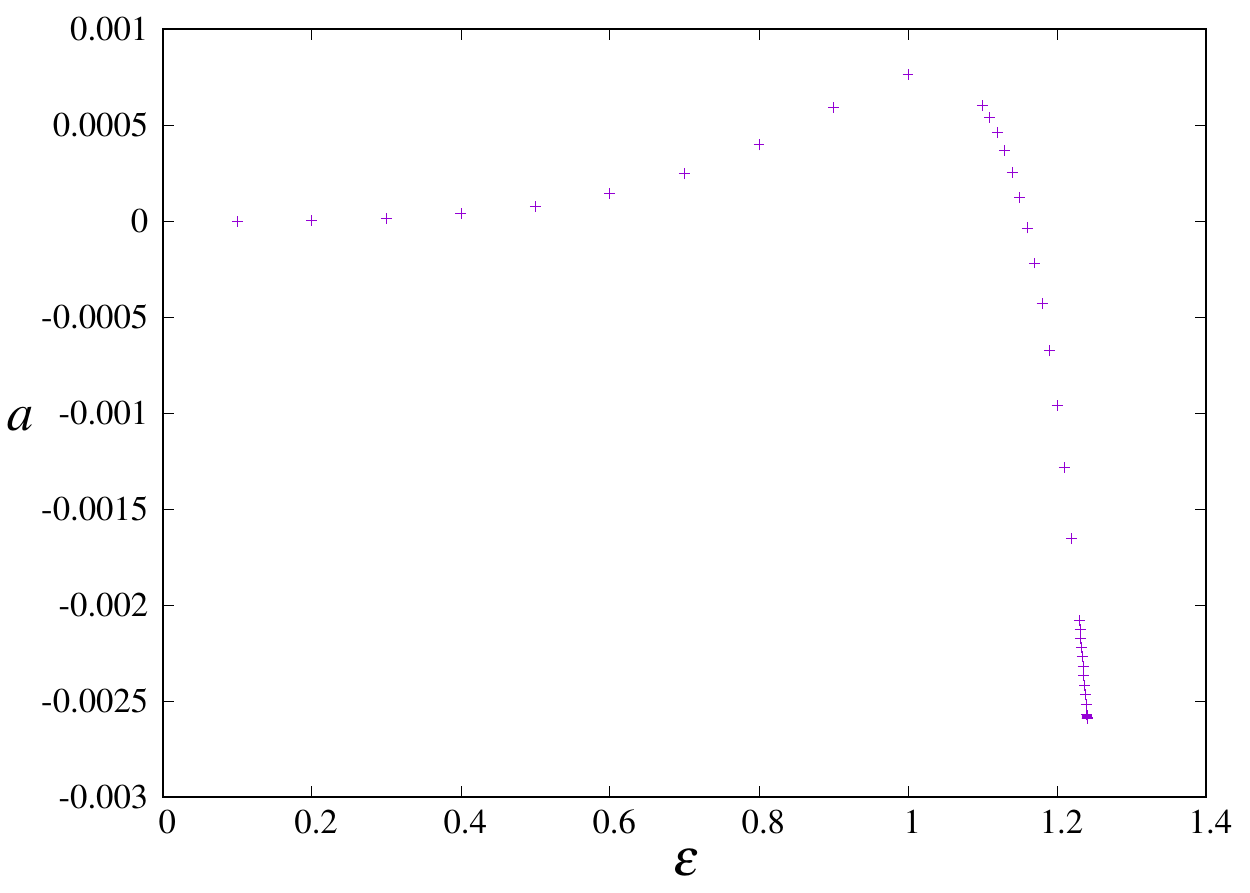} &
\includegraphics[width=0.35\linewidth]{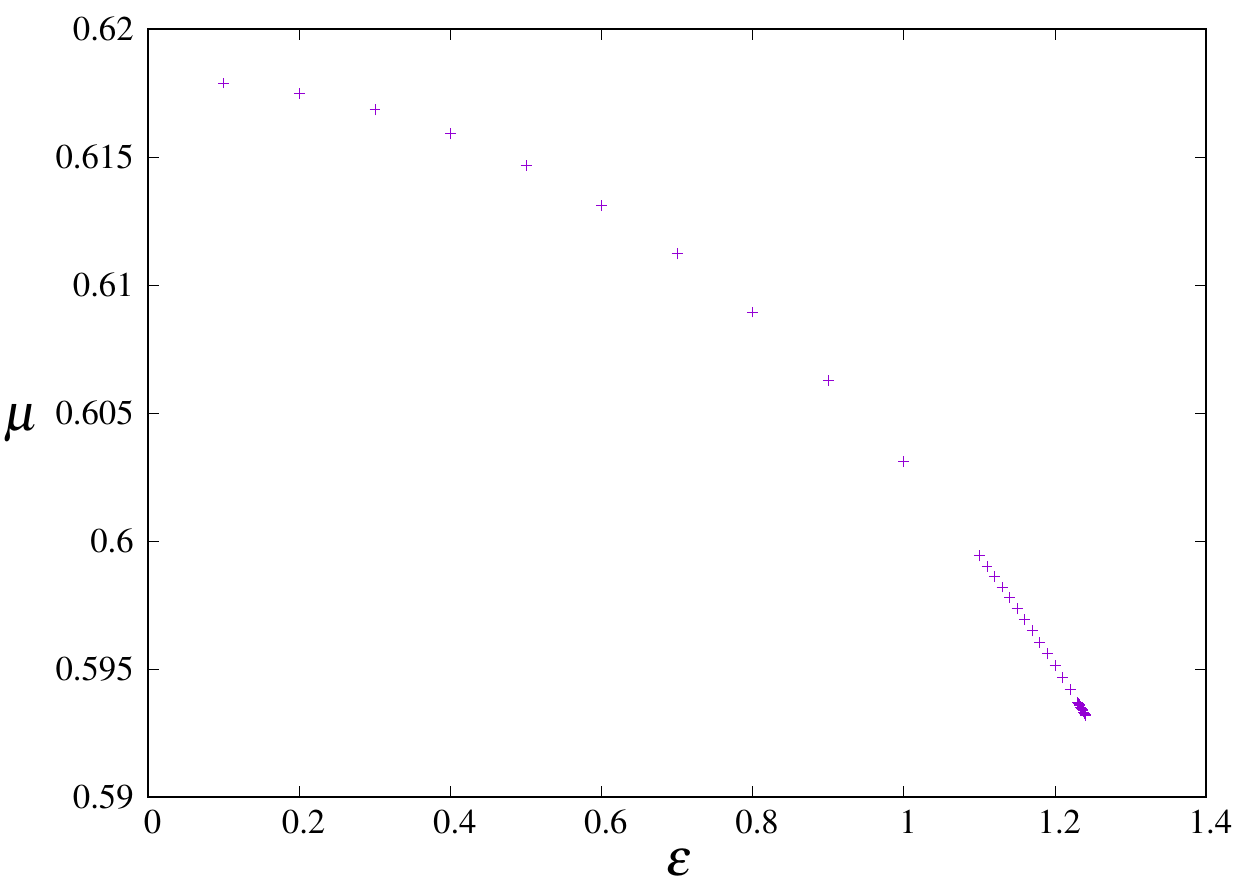}
\end{tabular}
\caption{\label{fig:ex2 parameters}
Continuation w.r.t. $\eps$ of a non-$a$-twist circle with frequency $\omega$ (non-symmetric case):
(left) adjusting parameter $a$; (left) unfolding parameter $\mu$.
}
\end{figure}

\begin{figure}[ht]
\begin{subfigure}[a][$\eps= 1.000000$, $a=  7.646104\cdot 10^{-4}$, $\mu= 0.6031124$]{ 
\centering
\begin{tabular}{cc}
\includegraphics[width=0.35\linewidth]{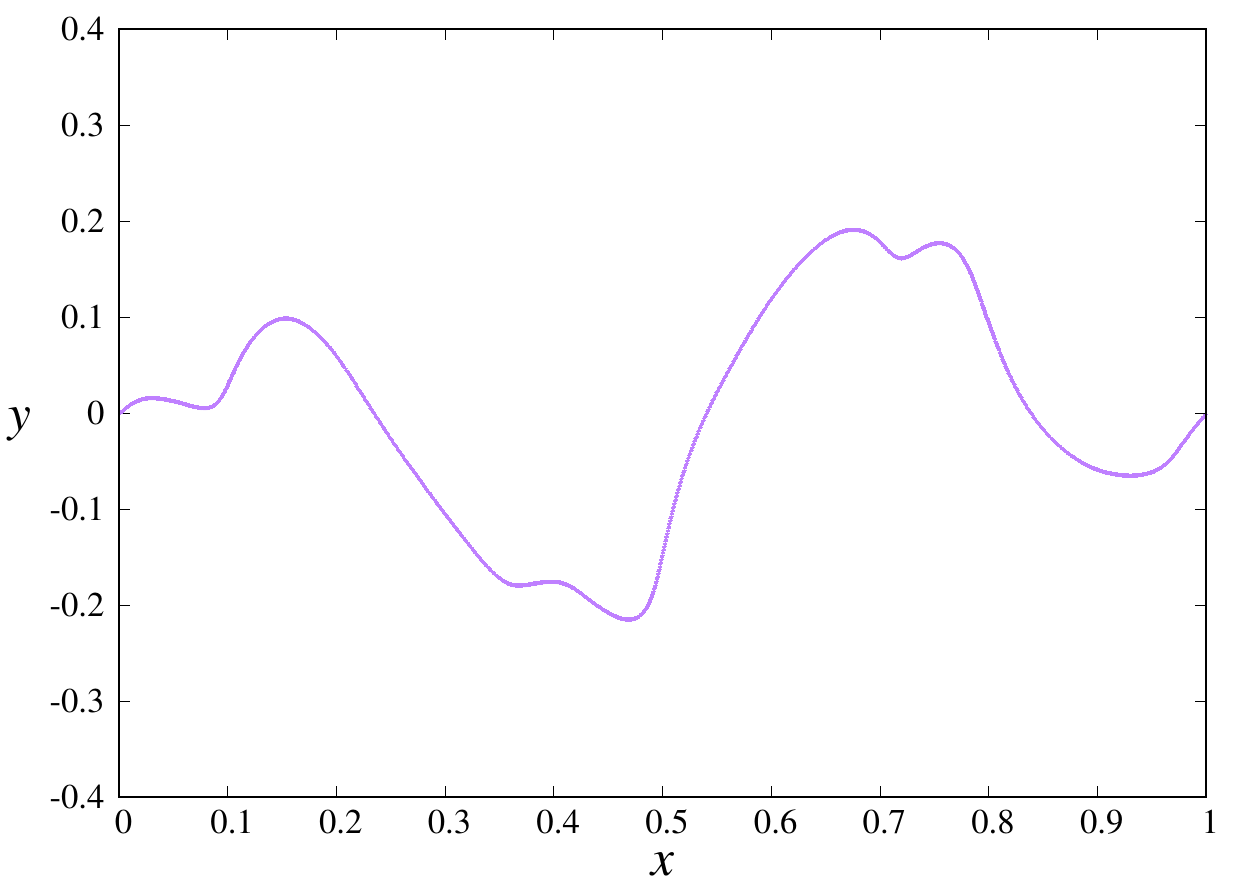} &
\includegraphics[width=0.35\linewidth]{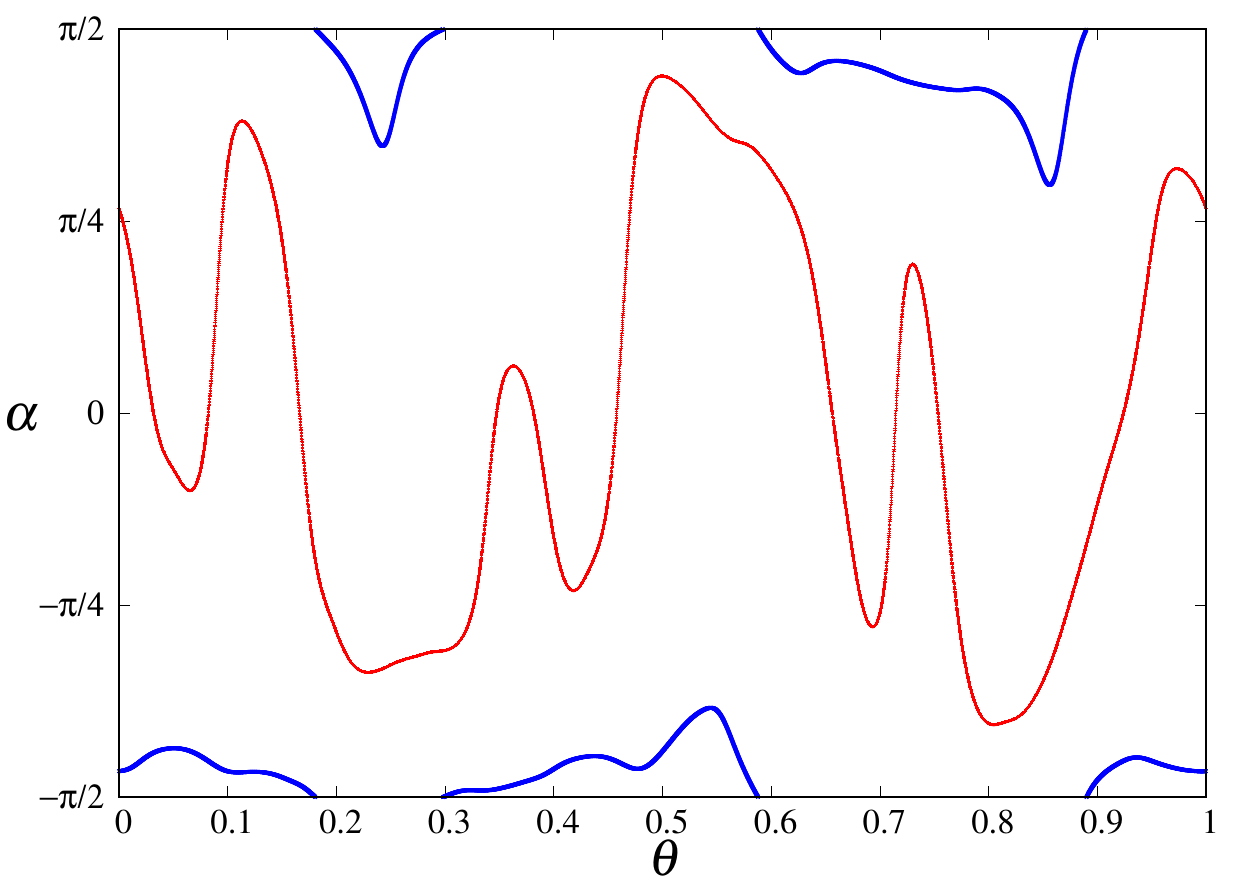} 
\end{tabular}
}
\end{subfigure}
\begin{subfigure}[a][$\eps= 1.200000$, $a=-9.571568\cdot 10^{-4}$, $\mu= 0.5951423$]{
\centering
\begin{tabular}{cc}
\includegraphics[width=0.35\linewidth]{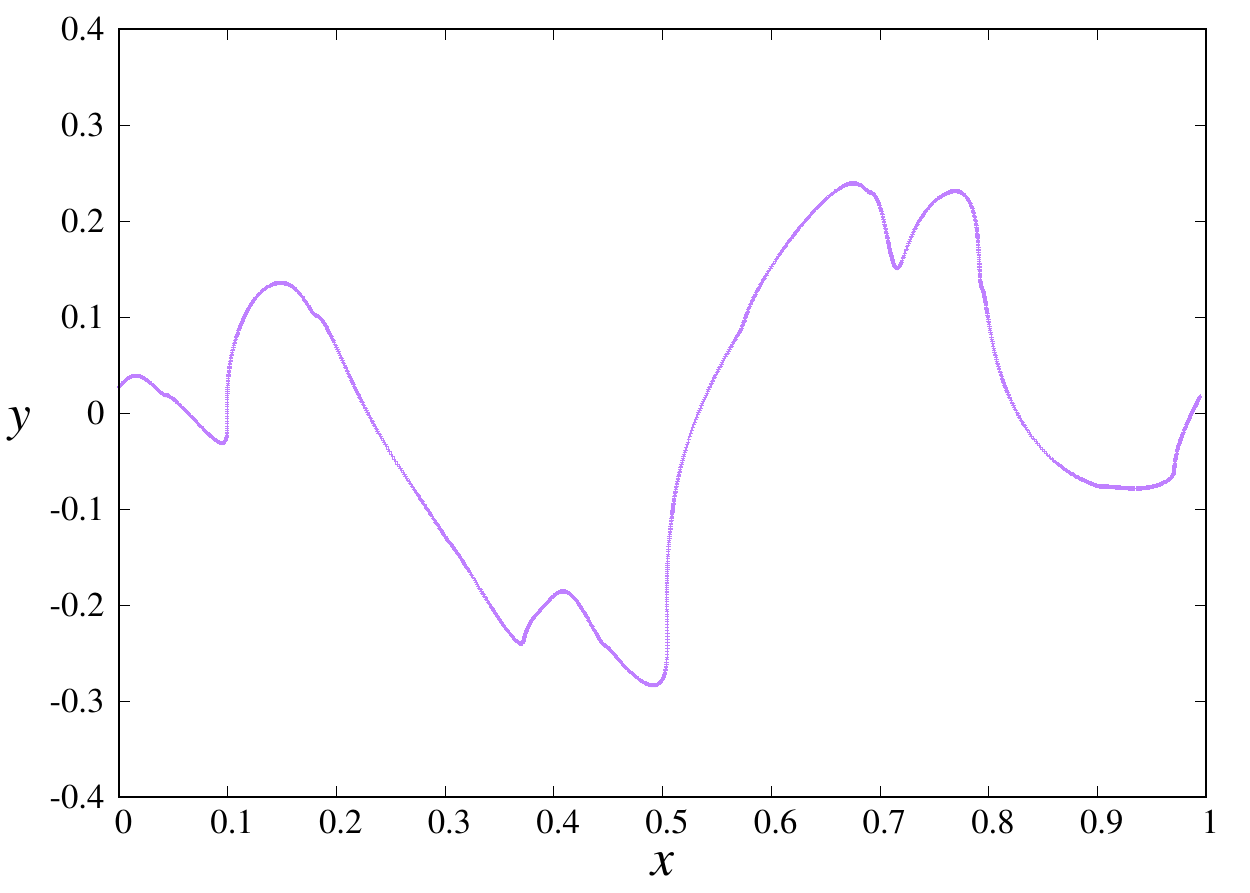} &
\includegraphics[width=0.35\linewidth]{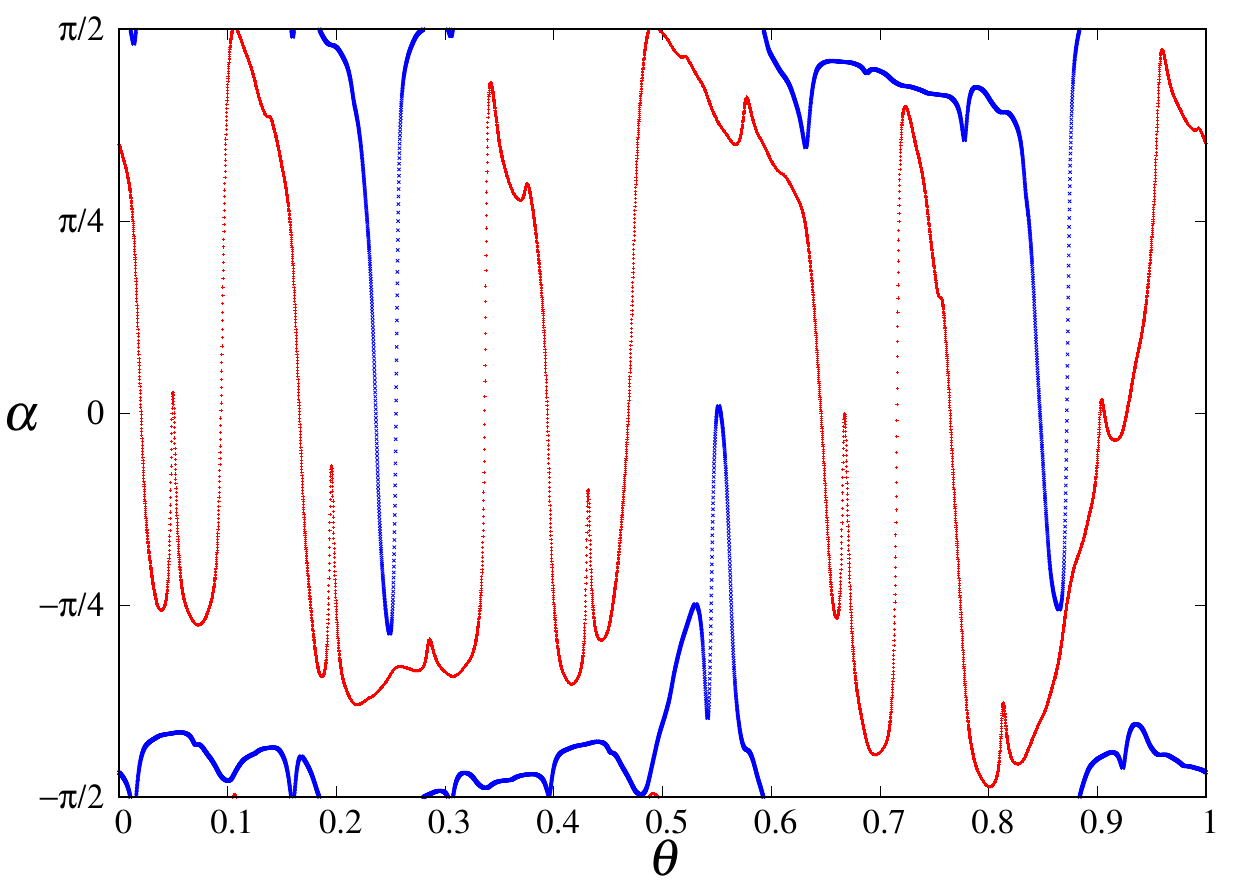} 
\end{tabular}
}
\end{subfigure}
\begin{subfigure}[a][$\eps= 1.240340$, $a=-2.588932\cdot 10^{-3}$, $\mu= 0.5932114$]{
\centering
\begin{tabular}{cc}
\includegraphics[width=0.35\linewidth]{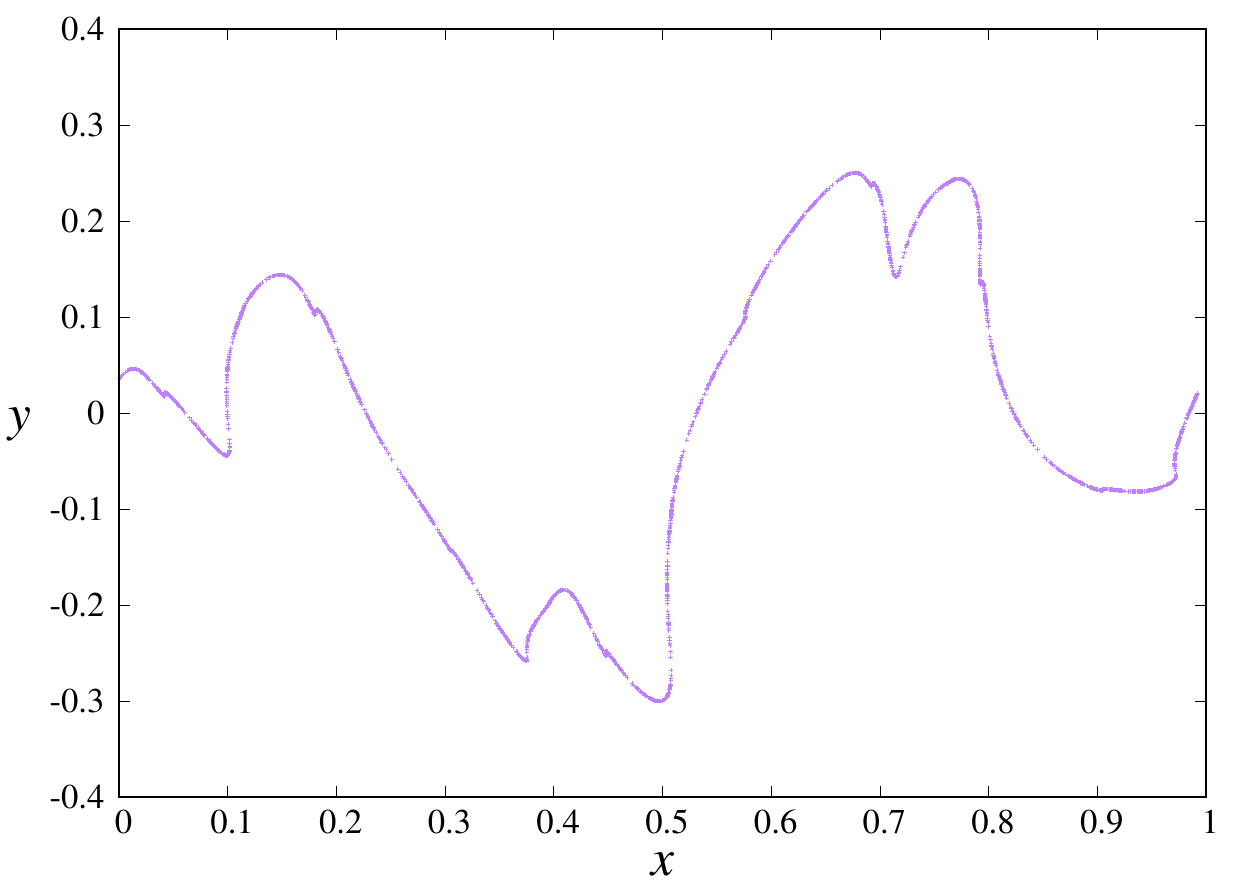} & 
\includegraphics[width=0.35\linewidth]{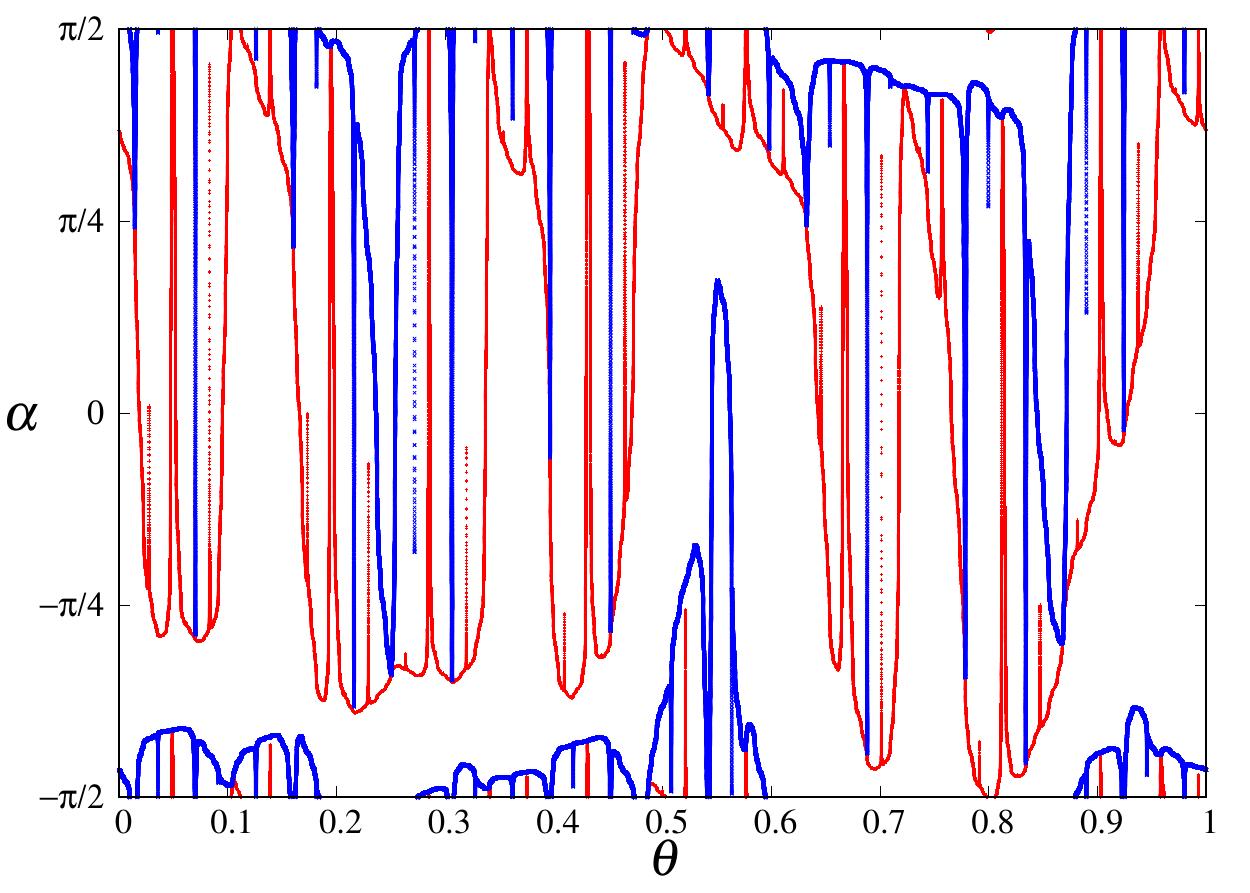}
\end{tabular}
}
\end{subfigure}
\caption{\label{fig:ex2 non-twist circle}
  Continuation w.r.t. $\eps$ of a non-$a$-twist circle with frequency $\omega$ (non-symmetric case):  
  (left) invariant circle; (right) projectivized tangent bundle (in red) and stable bundle (in blue).
  }
\end{figure}

\begin{figure}[ht]
\begin{tabular}{cc}
\includegraphics[width=0.35\linewidth]{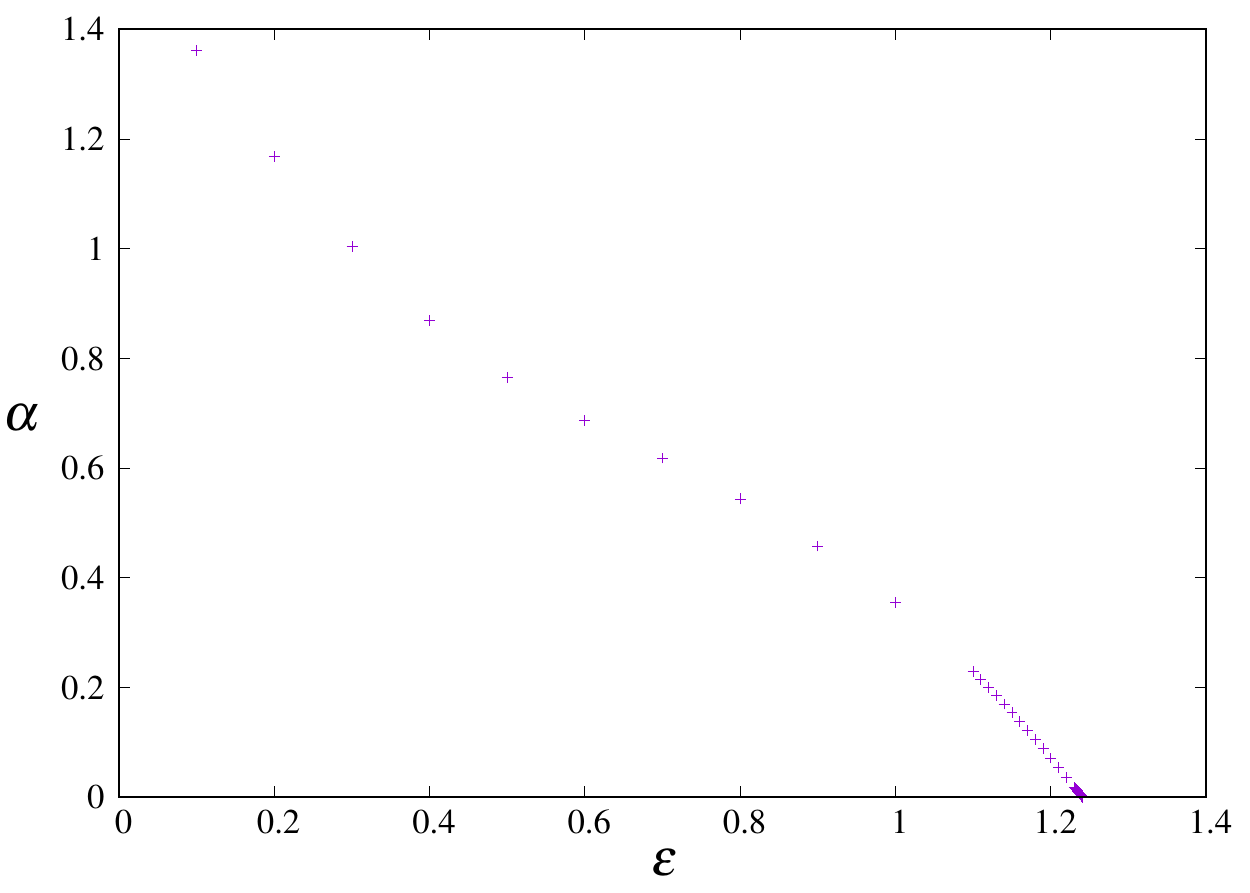} &
\includegraphics[width=0.35\linewidth]{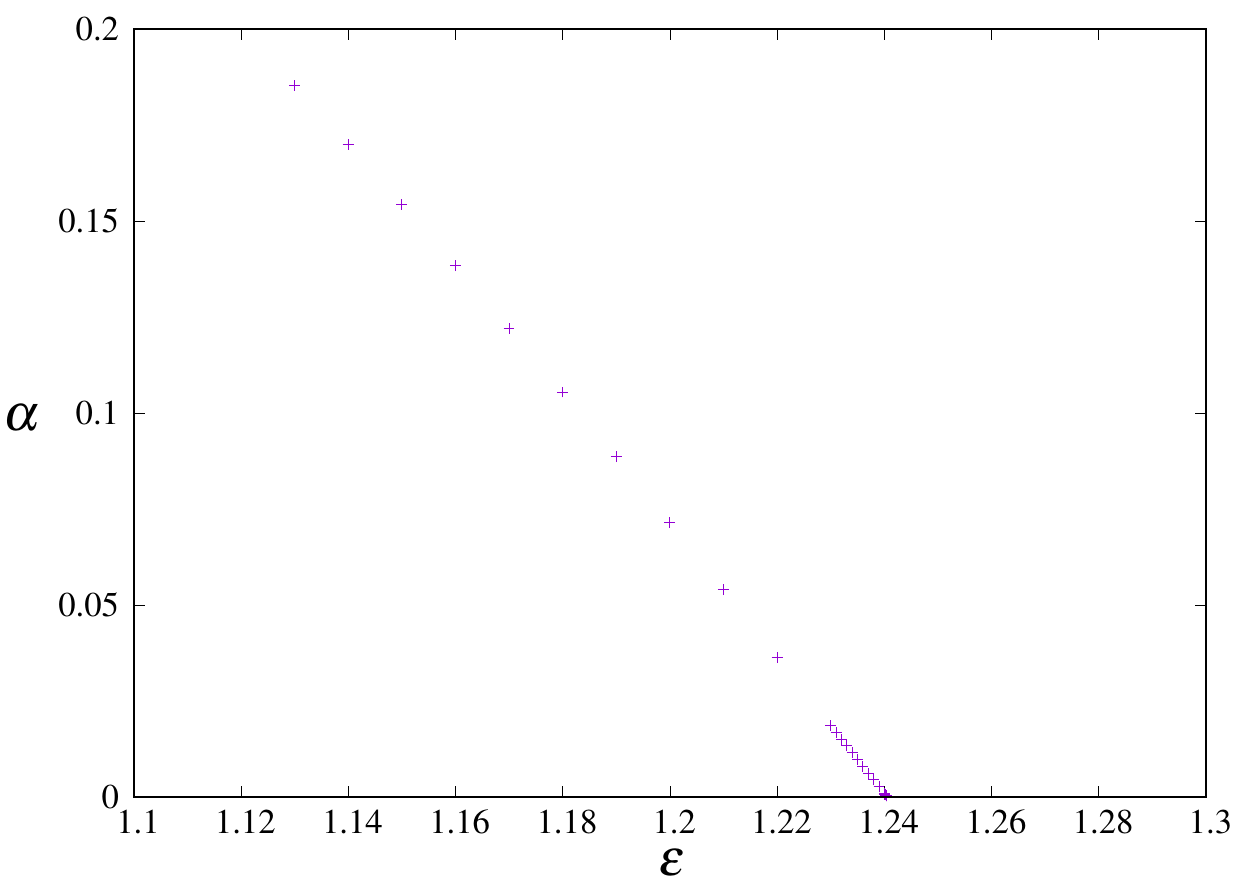}
\end{tabular}
\caption{\label{fig:ex2 angles}
Continuation w.r.t. $\eps$ of a non-$a$-twist circle with frequency $\omega$ (non-symmetric case): 
(left) minimum angle $\alpha$ between the stable and tangent bundles;
(right) critical behavior. The breakdown of the circle is produced at $\eps_{\rm c}\simeq 1.240522$.
}
\end{figure}

First, with algorithm of Section~\ref{sec:NT-tori}, we continue with respect to parameter $\eps$ a non-$a$-twist circle
with rotation number $\omega= \tfrac{1}{2} (\sqrt{5}-1)$. The adjusting parameters $a$ and $\mu$ as functions of
perturbation parameter $\eps$ are shown in Figure~\ref{fig:ex2 parameters}. Unlike the symmetric case, 
parameter $a$ varies, and remains bounded inside an interval of size $2.6\times 10^{-3}$ around
zero. The continuation reaches the value $\eps= 1.230340$, in which the invariant circle is 
approximated with a truncated Fourier series with $524288$ modes. The process of breakdown and the collision
of the invariant bundles is shown in Figure~\ref{fig:ex2 non-twist circle}.  
We notice that in contrast with the bundle collapse in the symmetric version
of the dissipative standard non-twist map, the collapse for this example
only happens on one side of the bundles, leaving a gap between the bundles.
The minimum angle between bundles is also asymptotically linear when close to breakdown, see 
Figure~\ref{fig:ex2 angles}, from which we can extrapolate the critical value $\eps_{\rm c}\simeq 1.240522$.

\begin{figure}
\begin{tabular}{cc}
\includegraphics[width=0.35\linewidth]{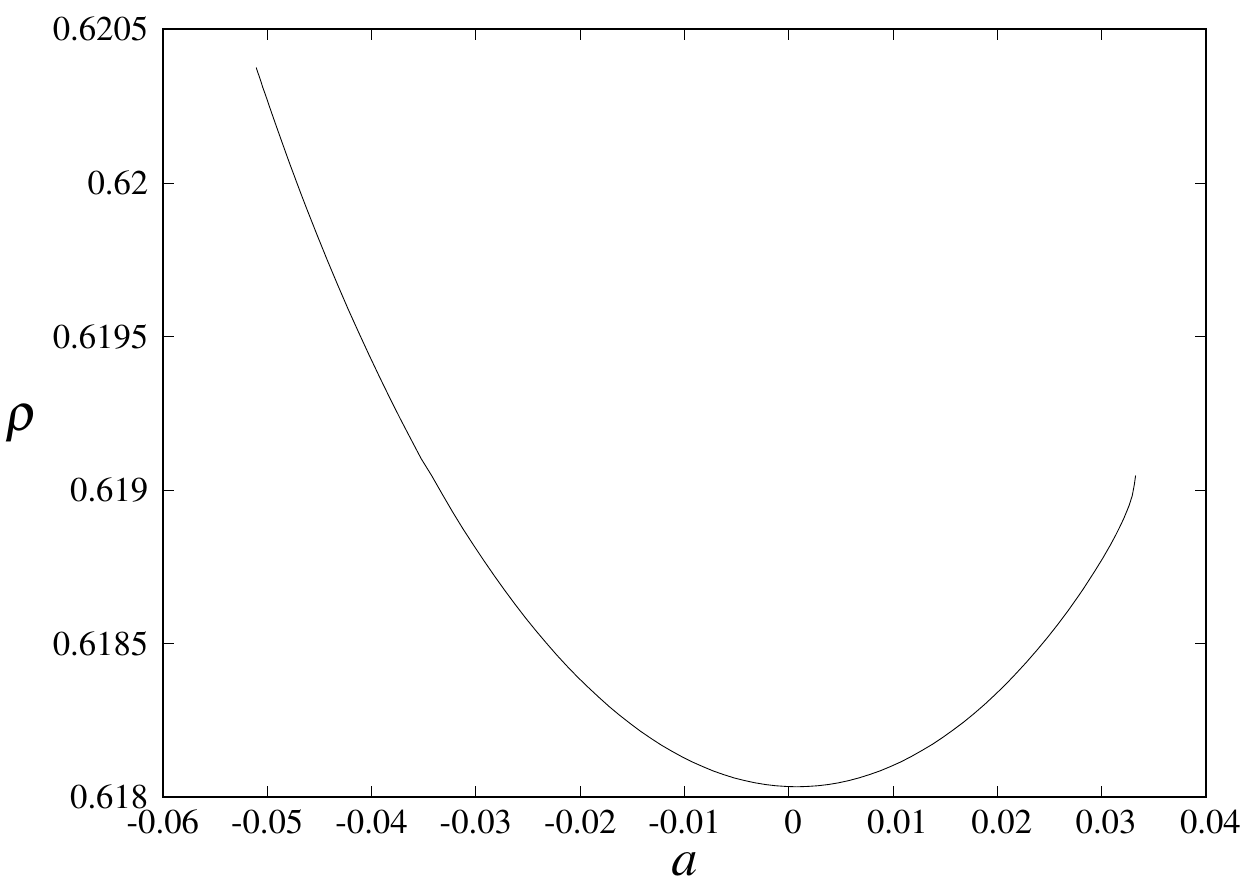} &
\includegraphics[width=0.35\linewidth]{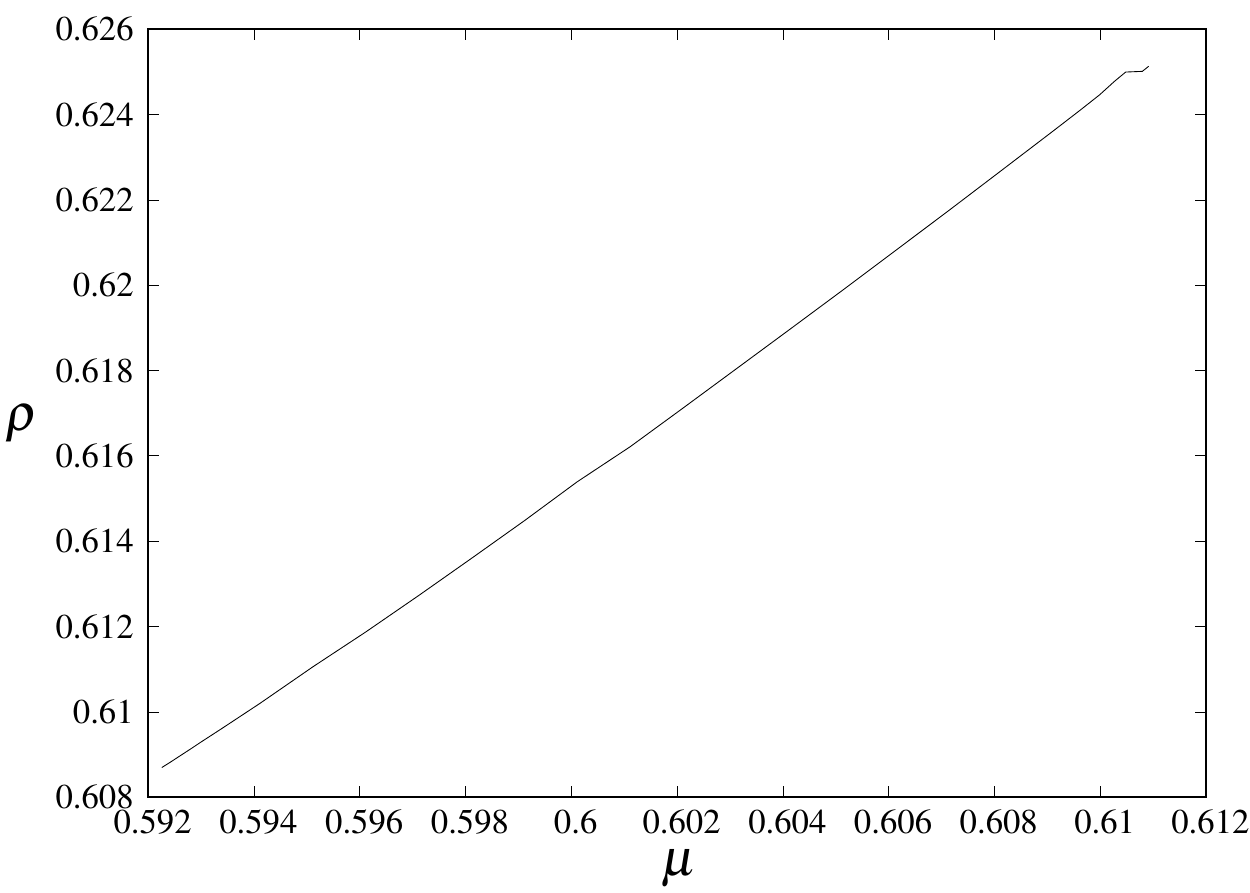}
\end{tabular}
\caption{\label{fig:ex2 rotation number}
Rotation number versus continuation parameter from the non-$a$-twist circle 
with $a= 7.646104\cdot 10^{-4}$, $\mu= 0.6031124$, $\eps= 1.00000$ (non-symmetric case):
(left) Continuation w.r.t. $a$; (right) Continuation w.r.t. $\mu$.
}
\end{figure}

As in the first example, in Figure~\ref{fig:ex2 rotation number} we show the graph of the rotation number 
of the invariant circle as a function of a parameter of  continuation (either $a$ or $\mu$) 
starting at a non-$a$-twist circle for $\eps= 1.00000$, $a= 7.646104\cdot 10^{-4}$, $\mu= 0.6031124$. 
The figure provides again a dynamical interpretation of the fact that the invariant
circle is non-$a$-twist, but $\mu$-twist.

\begin{figure}[ht]
  \includegraphics[width=0.75\linewidth]{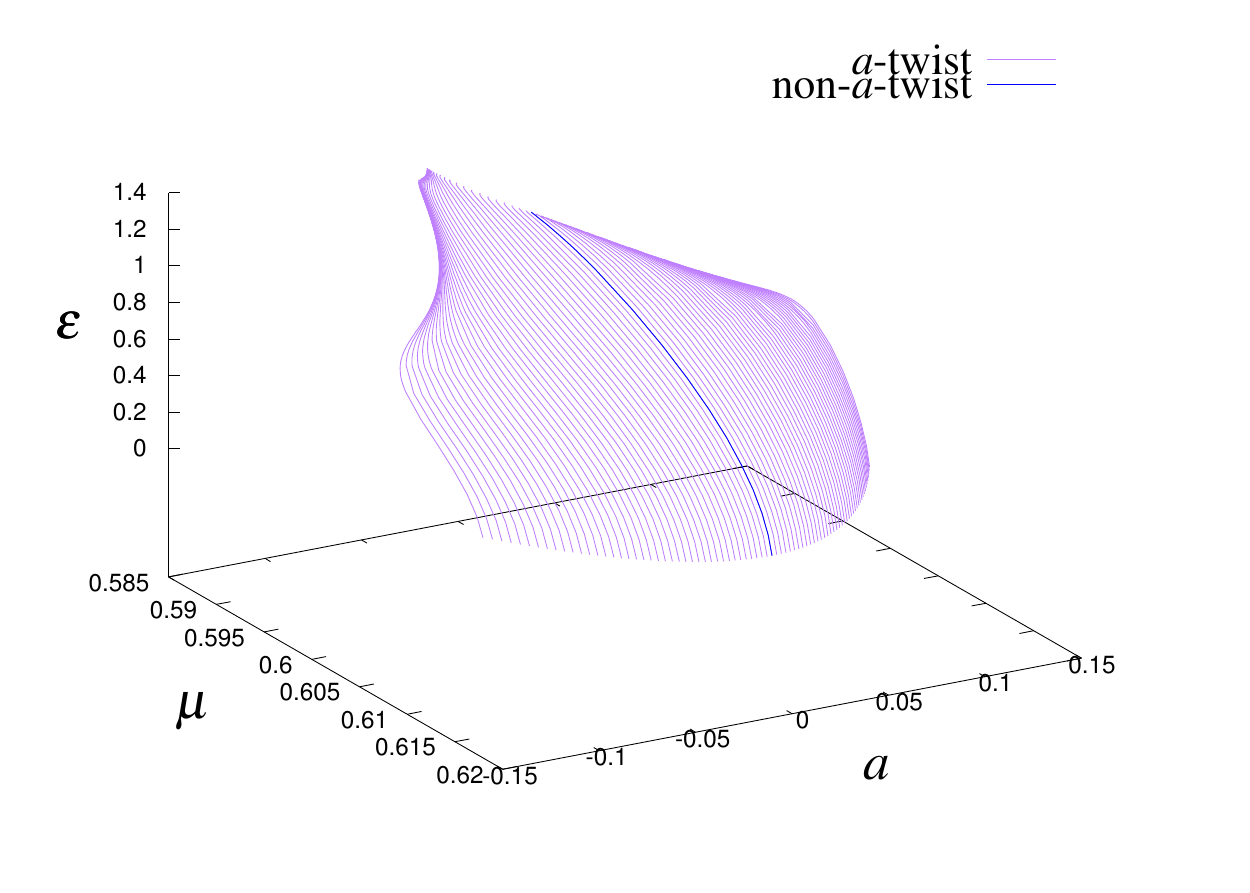}
\caption{\label{fig:ex2 surface}
Parameter surface generated by continuation w.r.t. $\eps$ of invariant circles with frequency $\omega$ and 
fix $a$-twist (non-symmetric case).
}
\end{figure}

In Figure~\ref{fig:ex2 surface}, we show continuations with 
respect to $\eps$ of invariant circles with fixed frequency $\omega$ and different
values of the $a$-twist. That is, we compute the surface of 
parameter points for which
there is an invariant circle with frequency $\omega$.
We plotted this surface showing the values of
$a$ and $\mu$ along the $\eps$ continuation. In particular, the continuation curve
  corresponding to an $a$-twist $b_a$ starts
  with $a= \frac12 b_a$, $\mu= \omega-a^2$ and $\eps= 0$. We have highligted
the curve corresponding
to zero $a$-twist. Note that the surface is not symmetric
with respect to $a$ and for negative $a$-twist there is a region where the
circles seem to persist for larger values of $\mu$ and $\eps$.

\section{Conclusions}

In this paper we have clarified  the property of being non-twist for a circle, 
in the context of conformally symplectic systems. 
This non-twist property has to do with the degeneracy condition arising when tuning a particular parameter to fix the dynamics of an invariant circle
to a given rotation number. Hence, the non-twist condition is with respect to a particular parameter.
As such, the concept can be extended to many other systems in which 
parameters have to be adjusted to fix the frequency, as in \cite{CanadellH17b}. In symplectic systems, the parameters
to adjust are the actions of a torus. 

We have also presented several algorithms for computing invariant circles, including non-twist circles
and a methodology to compute parametric surfaces in parameter space corresponding to invariant circles with a prescribed (Diophantine) frequency. The key of our methodology is introducing a concept of twist with respect to a parameter, 
so one can compute continuation curves corresponding to a fix twist. Unlike the symplectic case, 
non-twist tori in conformally symplectic
systems do not seem to be the more robust, meaning they are not the ones that survive for greater values of
perturbation parameters.

The algorithms are very efficient, and let us compute invariant circles even with hundreds of thousands of Fourier coefficients,
and then explore the regimes at the verge of analyticity breakdown. 

\section*{Acknowledgments}

R.C. was  partially supported by DGAPA-UNAM projects PAPIIT
IA 102818, IN101020 and by UIU project UCM-04-2019. M.C. was supported by MDM-2014-0445 (MINECO).
A.H. was supported by the grants PGC2018-100699-B-I00 (MCIU-AEI-FEDER, UE),
2017 SGR 1374 (AGAUR), MSCA 734557 (EU Horizon 2020), and MDM-2014-0445 (MINECO).
NSF under Grant No.~1440140  supported  R.C. and A.H., for their residences at 
MSRI in Berkeley, California, during the Fall 2018 semester.



\newcommand{\etalchar}[1]{$^{#1}$}
\def\cprime{$'$} \def\cprime{$'$} \def\cprime{$'$} \def\cprime{$'$}
\providecommand{\bysame}{\leavevmode\hbox to3em{\hrulefill}\thinspace}
\providecommand{\MR}{\relax\ifhmode\unskip\space\fi MR }
\providecommand{\MRhref}[2]{%
  \href{http://www.ams.org/mathscinet-getitem?mr=#1}{#2}
}
\providecommand{\href}[2]{#2}

\end{document}